\documentclass[11pt,a4paper]{article}
\usepackage{geometry}
\usepackage{amsmath,amssymb,amsthm,mathtools,bm}
\usepackage{thmtools}
\usepackage{booktabs,microtype,xcolor,fancyhdr}
\usepackage{graphicx,placeins,needspace,float}
\usepackage[font=small,labelfont=bf,hypcap=false]{caption}
\usepackage[ruled, noline, nofillcomment, linesnumbered]{algorithm2e}
\usepackage[hidelinks, hypertexnames=false]{hyperref}
\usepackage{bookmark}
\usepackage{cleveref}
\usepackage[square, comma, sort&compress, numbers]{natbib}

\hypersetup{pdftitle={Convergence Analysis of TriCG and TriMR for
Symmetric Quasi-Definite Systems}}
\numberwithin{equation}{section}

\DontPrintSemicolon
\SetAlgoCaptionSeparator{.}
\SetKwComment{tcc}{$\triangleright$\ }{}

\SetCommentSty{mycommfont}
\SetNlSty{textrm}{}{:}

\declaretheorem[
    name=Theorem,
    numberwithin=section,
    refname={theorem,theorems},
    Refname={Theorem,Theorems}
]{theorem}

\DeclareMathOperator{\range}{range}
\DeclareMathOperator{\rank}{rank}
\DeclareMathOperator{\diag}{diag}
\DeclareMathOperator{\spn}{span}
\DeclareMathOperator{\Null}{null}
\DeclareMathOperator*{\argmin}{{\arg\min}}
\newcommand{\R}{\mathbb{R}}
\newcommand{\Kry}{\mathcal{K}}
\newcommand{\one}{\mathbf{1}}
\newcommand{\TCG}{\mathrm{TriCG}}
\newcommand{\TMR}{\mathrm{TriMR}}
\newcommand{\BCG}{\mathrm{BCG}}
\newcommand{\BMR}{\mathrm{BMR}}

\title{Convergence analysis of TriCG and TriMR for symmetric
quasi-definite linear systems}

\author{Kui Du\thanks{School of Mathematical Sciences, Xiamen University, Xiamen
    361005, China (kuidu@xmu.edu.cn).}, \quad Jia-Jun Fan\thanks{School of
        Mathematical Sciences, Xiamen University, Xiamen 361005, China
(jiajunfan@stu.xmu.edu.cn).}}

\date{}

\begin{document}
\maketitle

\begin{abstract}
    TriCG and TriMR are two short-recurrence iterative methods introduced
    by Montoison and Orban [SIAM J.\ Sci.\ Comput., 43 (2021),
    pp.\ A2502--A2525] for solving symmetric quasi-definite
    linear systems. Extensive numerical experiments have demonstrated the
    efficiency of these two methods, but a systematic convergence analysis
    is still lacking. In this paper, we derive convergence bounds for the
    TriCG error norm and the TriMR residual norm.
    Numerical experiments support the theoretical results.
\end{abstract}



\section{Introduction}
\label{sec:introduction}
We consider the symmetric quasi-definite (SQD) linear system
\begin{equation}
    \label{eq:general-system}
    \begin{bmatrix}
        \mathbf{M} & \mathbf{A} \\
        \mathbf{A}^{\top} & - \mathbf{N}
    \end{bmatrix}
    \begin{bmatrix}
        \mathbf{x} \\
        \mathbf{y}
    \end{bmatrix}
    =
    \begin{bmatrix}
        \mathbf{b} \\
        \mathbf{c}
    \end{bmatrix},
\end{equation}
where $\mathbf{M} \in \R^{m \times m}$ and $\mathbf{N} \in \R^{n \times n}$
are symmetric positive definite matrices, and $\mathbf{A} \in \R^{m \times n}$
is a nonzero matrix, and
$\mathbf{b} \in \R^{m}$ and $\mathbf{c} \in \R^{n}$ are vectors.
The coefficient matrix is
indefinite and nonsingular even when $\mathbf{A}$ is rectangular or
rank deficient~\cite{vanderbei1995symmetric}. SQD linear systems
arise in a wide range of applications,
including optimization
problems~\cite{altman1999regularized,friedlander2012primal}
and saddle-point problems~\cite{benzi2005numerical}.

Montoison and Orban~\cite{montoison2021tricg} introduced TriCG and TriMR
for solving the SQD system \eqref{eq:general-system} with nonzero
$\mathbf{b}$ and $\mathbf{c}$, based on the generalized
Saunders--Simon--Yip (SSY) tridiagonalization~\cite{saunders1988two}.
TriCG and TriMR use Galerkin and residual-minimization principles,
respectively, analogous to those underlying the conjugate gradient
(CG)~\cite{hestenes1952methods} and minimal residual
(MINRES)~\cite{paige1975solutions} methods.
With a zero initial guess, each TriCG or TriMR iterate is the sum
of the two columns of the corresponding iterate obtained by applying
block-CG or block-MINRES to
\[
    \begin{bmatrix}
        \mathbf{M} & \mathbf{A} \\
        \mathbf{A}^{\top} & - \mathbf{N}
    \end{bmatrix}
    \begin{bmatrix}
        \mathbf{x}^{b} & \mathbf{x}^{c} \\
        \mathbf{y}^{b} & \mathbf{y}^{c}
    \end{bmatrix}
    =
    \begin{bmatrix}
        \mathbf{b} & \mathbf{0} \\
        \mathbf{0} & \mathbf{c}
    \end{bmatrix}
\]
with the preconditioner
$\diag(\mathbf{M},\mathbf{N})$~\cite{montoison2021tricg}.
For the special case in which $\mathbf{b} = \mathbf{0}$ or
$\mathbf{c} = \mathbf{0}$,
Orban and Arioli~\cite{orban2017iterative} developed generalizations of
LSQR~\cite{paige1982lsqr}, LSMR~\cite{fong2011lsmr}, and Craig's
method~\cite{craig1955n} using the generalized Golub--Kahan
bidiagonalization process~\cite{arioli2013generalized}.

Du, Fan, and Zhang \cite{du2025improved} recently established the
finite termination of the generalized SSY process, with the maximum
number of iterations determined by the
rank of $\mathbf{A}$ and the number of its distinct elliptic singular
values. They also demonstrated that TriCG and TriMR may fail to
produce the solution when the generalized SSY process terminates, and
proposed improved versions of TriCG and TriMR to address this issue. Numerical
experiments indicated that TriCG can converge
slowly when $\mathbf{A}$ has a substantial number of large elliptic
singular values. Du and Fan~\cite{du2026deflated} developed a variant
of TriCG with deflated restarting to accelerate convergence. Despite
these advances, convergence bounds for TriCG and TriMR have not yet
been established, to the best of our knowledge. In this paper, we
derive bounds for the TriCG error norm and the TriMR residual norm.

Although the TriCG and TriMR iterates are related to those of
block-CG and block-MINRES, the convergence bounds obtained in
\cite{oleary1980block,simoncini1996convergence} rely on
spectral assumptions that are not satisfied by SQD linear systems.
Consequently, the existing convergence bounds for block-CG and
block-MINRES cannot be used to explain the practical performance
of TriCG and TriMR on SQD linear systems. In this paper, we fill this
gap. Specifically, we first analyze the special case
$\mathbf{M} = \mathbf{I}_{m}$ and $\mathbf{N} = \mathbf{I}_{n}$.
We show that the TriCG error norm is bounded by the TriMR residual
norm at each iteration. Using this relation, we derive convergence
bounds for both methods: a bound analogous to the classical CG bound
and two refinements that exploit the block structure and the singular
values of the off-diagonal block to obtain potentially smaller
convergence factors. We then extend these results to general SQD linear
systems in the corresponding weighted norms. Numerical experiments
compare the derived bounds with the observed convergence behavior of
TriCG and TriMR.

This paper is organized as follows. The remainder of this section
introduces the notation used throughout the paper.
\Cref{sec:background} reviews the generalized SSY process, and
\Cref{subsec:tricg-trimr} describes TriCG and TriMR.
\Cref{sec:common} presents some results needed for the convergence
analysis of TriCG and TriMR. \Cref{sec:even,sec:odd} derive
convergence bounds for even and odd iterations, respectively.
\Cref{sec:extension} extends these bounds to general SQD linear
systems in weighted norms. Numerical experiments and conclusions are
given in \Cref{sec:numerical} and \Cref{sec:discussion}, respectively.

\emph{Notation.}
Matrices and vectors are denoted by bold uppercase and lowercase
letters, respectively; scalars are italic.
The symbols $\mathbf{I}_{d}$, $\mathbf{e}_{i}$, and $\one_{d}$ denote the
identity matrix of order $d$, the $i$th standard basis vector, and
the all-ones vector of length $d$. The symbol $\mathbf{0}$ denotes
a zero vector or matrix of the appropriate size. Dimensions are
omitted when they are clear from the context.
For a matrix $\mathbf{A}$, $\mathbf{A}^{\top}$ denotes its transpose,
$\mathbf{A}^{-1}$ its inverse when it is nonsingular, and
$\mathbf{A}^{\dagger}$ its Moore--Penrose pseudoinverse.
Unless otherwise specified, $\lVert \cdot \rVert$ denotes the Euclidean
vector norm or the induced matrix $2$-norm. For a symmetric positive
definite matrix $\mathbf{M}$,
$\mathbf{M}^{1/2}$ denotes its unique symmetric positive definite square
root, satisfying $(\mathbf{M}^{1/2})^2=\mathbf{M}$, and
$\mathbf{M}^{-1/2}=(\mathbf{M}^{1/2})^{-1}$.
The associated weighted norm is
$\lVert
\mathbf{v} \rVert_{\mathbf{M}} = (\mathbf{v}^{\top} \mathbf{M}
\mathbf{v})^{1/2}$.
The notation $\range(\mathbf{A})$, $\Null(\mathbf{A})$, and
$\rank(\mathbf{A})$ refers to the column space, nullspace, and rank
of $\mathbf{A}$. The span of a collection of vectors, or of the
columns of a collection of matrices, is denoted by $\spn \{ \cdot \}$.
The operator $\diag(\cdot)$ forms a diagonal or block-diagonal
matrix from its arguments.

\section{The generalized Saunders--Simon--Yip process}
\label{sec:background}
\label{subsec:gssy}
For symmetric positive definite matrices $\mathbf{M} \in \R^{m \times m}$
and $\mathbf{N} \in \R^{n \times n}$, a matrix
$\mathbf{A} \in \R^{m \times n}$, and nonzero vectors
$\mathbf{b} \in \R^{m}$ and $\mathbf{c} \in \R^{n}$,
\Cref{alg:orth_tri} summarizes the generalized Saunders--Simon--Yip
(SSY) tridiagonalization process~\cite{saunders1988two,montoison2021tricg}.

\begin{algorithm}[htb]
    \caption{Generalized Saunders--Simon--Yip tridiagonalization process}
    \label{alg:orth_tri}
    \KwIn{Symmetric positive definite matrices
        $\mathbf{M} \in \R^{m \times m}$ and $\mathbf{N} \in \R^{n \times n}$,
        a matrix $\mathbf{A} \in \R^{m \times n}$, and nonzero vectors
    $\mathbf{b} \in \R^{m}$ and $\mathbf{c} \in \R^{n}$}
    $\mathbf{u}_{0} = \mathbf{0}$, $\mathbf{v}_{0} = \mathbf{0}$\;
    $\widetilde{\mathbf{u}}_{1} = \mathbf{M}^{-1} \mathbf{b}$,
    $\widetilde{\mathbf{v}}_{1} = \mathbf{N}^{-1} \mathbf{c}$\;
    $\beta_{1} = \sqrt{\widetilde{\mathbf{u}}_{1}^{\top} \mathbf{b}}$,
    $\gamma_{1} = \sqrt{\widetilde{\mathbf{v}}_{1}^{\top} \mathbf{c}}$\;
    $\mathbf{u}_{1} = \widetilde{\mathbf{u}}_{1}/ \beta_{1}$,
    $\mathbf{v}_{1} = \widetilde{\mathbf{v}}_{1}/ \gamma_{1}$\;
    \For{$k = 1, 2, \dots$}{
        $\mathbf{q} = \mathbf{A} \mathbf{v}_{k} - \gamma_{k} \mathbf{M}
        \mathbf{u}_{k - 1}$\;
        $\mathbf{p} = \mathbf{A}^{\top} \mathbf{u}_{k} - \beta_{k}
        \mathbf{N} \mathbf{v}_{k - 1}$\;
        $\alpha_{k} = \mathbf{u}_{k}^{\top} \mathbf{q}$\;
        $\mathbf{q} = \mathbf{q} - \alpha_{k} \mathbf{M} \mathbf{u}_{k}$,
        $\mathbf{p} = \mathbf{p} - \alpha_{k} \mathbf{N} \mathbf{v}_{k}$\;
        $\widetilde{\mathbf{u}}_{k + 1} = \mathbf{M}^{-1} \mathbf{q}$,
        $\widetilde{\mathbf{v}}_{k + 1} = \mathbf{N}^{-1} \mathbf{p}$\;
        $\beta_{k + 1} = \sqrt{\widetilde{\mathbf{u}}_{k + 1}^{\top}
        \mathbf{q}}$,
        $\gamma_{k + 1} = \sqrt{\widetilde{\mathbf{v}}_{k + 1}^{\top}
        \mathbf{p}}$\;
        \lIf{$\beta_{k + 1} = 0$ or $\gamma_{k + 1} = 0$}{\textbf{stop}}
        $\mathbf{u}_{k + 1} = \widetilde{\mathbf{u}}_{k + 1}/ \beta_{k + 1}$,
        $\mathbf{v}_{k + 1} = \widetilde{\mathbf{v}}_{k + 1}/ \gamma_{k + 1}$\;
    }
\end{algorithm}

Define the basis matrices
\[
    \mathbf{U}_{k} =
    \begin{bmatrix}
        \mathbf{u}_{1} & \cdots & \mathbf{u}_{k}
    \end{bmatrix} \in \R^{m \times k}, \quad
    \mathbf{V}_{k} =
    \begin{bmatrix}
        \mathbf{v}_{1} & \cdots & \mathbf{v}_{k}
    \end{bmatrix} \in \R^{n \times k},
\]
and
\[
    \mathbf{T}_{k} =
    \begin{bmatrix}
        \alpha_{1} & \gamma_{2} \\
        \beta_{2} & \alpha_{2} & \ddots \\
        & \ddots & \ddots & \gamma_{k} \\
        & & \beta_{k} & \alpha_{k}
    \end{bmatrix} \in \R^{k \times k}, \quad
    \mathbf{T}_{k + 1, k} =
    \begin{bmatrix}
        \mathbf{T}_{k} \\
        \beta_{k + 1} \mathbf{e}_{k}^{\top}
    \end{bmatrix}, \quad
    \mathbf{T}_{k, k + 1} =
    \begin{bmatrix}
        \mathbf{T}_{k} & \gamma_{k + 1} \mathbf{e}_{k}
    \end{bmatrix}.
\]
In exact arithmetic, assuming that no termination has occurred after
$k$ iterations, \Cref{alg:orth_tri} gives the following relations:
\begin{equation}\label{eq:gssy}
    \begin{gathered}
        \mathbf{A} \mathbf{V}_{k}
        = \mathbf{M} \mathbf{U}_{k} \mathbf{T}_{k}
        + \beta_{k + 1} \mathbf{M} \mathbf{u}_{k + 1} \mathbf{e}_{k}^{\top}
        = \mathbf{M} \mathbf{U}_{k + 1} \mathbf{T}_{k + 1, k}, \\
        \mathbf{A}^{\top} \mathbf{U}_{k}
        = \mathbf{N} \mathbf{V}_{k} \mathbf{T}_{k}^{\top}
        + \gamma_{k + 1} \mathbf{N} \mathbf{v}_{k + 1} \mathbf{e}_{k}^{\top}
        = \mathbf{N} \mathbf{V}_{k + 1} \mathbf{T}_{k, k + 1}^{\top}, \\
        \mathbf{U}_{k}^{\top} \mathbf{M} \mathbf{U}_{k}
        = \mathbf{V}_{k}^{\top} \mathbf{N} \mathbf{V}_{k} = \mathbf{I}_{k},
        \quad \mathbf{U}_{k}^{\top} \mathbf{A} \mathbf{V}_{k} = \mathbf{T}_{k}.
    \end{gathered}
\end{equation}

\section{TriCG and TriMR}
\label{subsec:tricg-trimr}
We assume that all methods use the zero initial guess.
For brevity, define
\begin{equation*}
    \begin{gathered}
        \mathbf{K} =
        \begin{bmatrix}
            \mathbf{M} & \mathbf{A} \\
            \mathbf{A}^{\top} & - \mathbf{N}
        \end{bmatrix}, \quad
        \mathbf{H} =
        \begin{bmatrix}
            \mathbf{M} & \mathbf{0} \\
            \mathbf{0} & \mathbf{N}
        \end{bmatrix}, \quad
        \mathbf{B}_{0} =
        \begin{bmatrix}
            \mathbf{b} & \mathbf{0} \\
            \mathbf{0} & \mathbf{c}
        \end{bmatrix}, \quad
        \mathbf{s} =
        \begin{bmatrix}
            \mathbf{x} \\
            \mathbf{y}
        \end{bmatrix}, \quad
        \mathbf{f} =
        \begin{bmatrix}
            \mathbf{b} \\
            \mathbf{c}
        \end{bmatrix} \ne \mathbf{0}.
    \end{gathered}
\end{equation*}
Let $\mathbf{U}_{k}$ and $\mathbf{V}_{k}$ be the basis matrices generated by
\Cref{alg:orth_tri}, and define
\begin{equation}
    \label{eq:searchspace}
    \mathbf{W}_{k} =
    \begin{bmatrix}
        \mathbf{U}_{k} & \mathbf{0} \\
        \mathbf{0} & \mathbf{V}_{k}
    \end{bmatrix}.
\end{equation}

TriCG seeks the iterate $\mathbf{s}_{k}^{\TCG} = \mathbf{W}_{k}
\mathbf{z}_{k}^{\TCG}$ which satisfies the Galerkin condition
\begin{equation} \label{eq:tricg_galerkin}
    \mathbf{W}_{k}^{\top}(\mathbf{f} - \mathbf{K}
    \mathbf{s}_{k}^{\TCG}) = \mathbf{0}.
\end{equation}
By \eqref{eq:gssy}, this leads to the subproblem
\begin{equation*}
    \begin{bmatrix}
        \mathbf{I}_{k} & \mathbf{T}_{k} \\
        \mathbf{T}_{k}^{\top} & - \mathbf{I}_{k}
    \end{bmatrix} \mathbf{z}_{k}^{\TCG}
    =
    \begin{bmatrix}
        \beta_{1} \mathbf{e}_{1} \\
        \gamma_{1} \mathbf{e}_{1}
    \end{bmatrix}.
\end{equation*}
TriMR seeks the iterate $\mathbf{s}_{k}^{\TMR} = \mathbf{W}_{k}
\mathbf{z}_{k}^{\TMR}$
which minimizes the $\mathbf{H}^{-1}$-norm of the residual
\cite[Section 3.2]{montoison2021tricg}
\begin{equation} \label{eq:trimrsmall}
    \min_{\mathbf{z} \in \R^{2 k}}
    \lVert \mathbf{f} - \mathbf{K} \mathbf{W}_{k} \mathbf{z}
    \rVert_{\mathbf{H}^{-1}}.
\end{equation}
The corresponding subproblem is
\[
    \mathbf{z}_{k}^{\TMR} = \argmin_{\mathbf{z} \in \R^{2 k}}
    \biggl\lVert
    \begin{bmatrix}
        \mathbf{I}_{k + 1, k} & \mathbf{T}_{k + 1, k} \\
        \mathbf{T}_{k, k + 1}^{\top} & - \mathbf{I}_{k + 1, k}
    \end{bmatrix} \mathbf{z}
    -
    \begin{bmatrix}
        \beta_{1} \mathbf{e}_{1} \\
        \gamma_{1} \mathbf{e}_{1}
    \end{bmatrix}
    \biggr\rVert,
\]
where $\mathbf{I}_{k + 1, k}$ consists of the first $k$ columns of
$\mathbf{I}_{k + 1}$.

Montoison and Orban~\cite{montoison2021tricg} established the equivalence
of TriCG and TriMR with block-CG and block-MINRES, respectively.
In particular, the $k$th iterate
$\mathbf{X}_{k}^{\BCG} = \mathbf{W}_{k} \mathbf{Z}_{k}^{\BCG}$ of block-CG
applied to $\mathbf{K} \mathbf{X} = \mathbf{B}_{0}$ with preconditioner
$\mathbf{H}$ satisfies the Galerkin condition
\[
    \mathbf{W}_{k}^{\top}(\mathbf{B}_{0} - \mathbf{K} \mathbf{X}_{k}^{\BCG})
    = \mathbf{0}.
\]
The corresponding subproblem is
\[
    \begin{bmatrix}
        \mathbf{I}_{k} & \mathbf{T}_{k} \\
        \mathbf{T}_{k}^{\top} & - \mathbf{I}_{k}
    \end{bmatrix} \mathbf{Z}_{k}^{\BCG}
    =
    \begin{bmatrix}
        \beta_{1} \mathbf{e}_{1} & \mathbf{0} \\
        \mathbf{0} & \gamma_{1} \mathbf{e}_{1}
    \end{bmatrix}.
\]
The TriCG iterate is the sum of the two columns of the block-CG iterate:
\begin{equation*}
    \mathbf{s}_{k}^{\TCG} = \mathbf{X}_{k}^{\BCG} \one_{2}, \quad
    \one_{2} =
    \begin{bmatrix}
        1 & 1
    \end{bmatrix}^{\top}.
\end{equation*}
An analogous argument using the residual minimization condition gives
$\mathbf{s}_{k}^{\TMR} = \mathbf{X}_{k}^{\BMR} \one_{2}$, where
$\mathbf{X}_{k}^{\BMR}$ is the corresponding block-MINRES iterate.

In exact arithmetic, TriCG and TriMR inherit the finite termination
property of the generalized SSY process. However, TriCG and TriMR may
fail to yield the solution when the generalized SSY process
terminates. To address this issue, Du et al.\ \cite{du2025improved}
proposed improved versions of TriCG and TriMR.

Despite the relations between the iterates of TriCG and TriMR and those of
block-CG and block-MINRES, the convergence bounds considered in
\cite{oleary1980block,simoncini1996convergence}
do not apply directly to SQD linear systems. The block-CG
error bound in \cite[Theorem 5]{oleary1980block} assumes a symmetric positive
definite coefficient
matrix. For symmetric linear systems, block-MINRES is equivalent to
block-GMRES. The block-GMRES residual bound in \cite[Theorem
2.2]{simoncini1996convergence} assumes that all eigenvalues have
positive real parts.
Neither assumption is satisfied by the symmetric indefinite
coefficient matrix of an SQD linear system.

\section{Preliminary results for the convergence analysis}
\label{sec:common}
In this section and the subsequent convergence analysis, we consider
the special case $\mathbf{M} = \mathbf{I}_{m}$ and $\mathbf{N} =
\mathbf{I}_{n}$, for which \eqref{eq:general-system} becomes
\begin{equation}
    \label{eq:identity-system}
    \begin{bmatrix}
        \mathbf{I}_{m} & \mathbf{A} \\
        \mathbf{A}^{\top} & - \mathbf{I}_{n}
    \end{bmatrix}
    \begin{bmatrix}
        \mathbf{x} \\ \mathbf{y}
    \end{bmatrix} =
    \begin{bmatrix}
        \mathbf{b} \\ \mathbf{c}
    \end{bmatrix}.
\end{equation}
The results extend to general symmetric
quasi-definite systems, as explained in \Cref{sec:extension}.
Throughout the analysis, we work in exact arithmetic and assume that the
generalized SSY process generates the required basis matrices with full
column rank.

\subsection{Eigenvalue distribution}
\label{subsec:eigenvalue-distribution}
The singular value decomposition (SVD) of $\mathbf{A}$ gives
\begin{equation}
    \label{eq:svd}
    \mathbf{A} =
    \begin{bmatrix}
        \widehat{\mathbf{U}}_{r} & \widehat{\mathbf{U}}_{c}
    \end{bmatrix}
    \begin{bmatrix}
        \mathbf{\Sigma}_{r} & \mathbf{0} \\
        \mathbf{0} & \mathbf{0}
    \end{bmatrix}
    \begin{bmatrix}
        \widehat{\mathbf{V}}_{r} & \widehat{\mathbf{V}}_{c}
    \end{bmatrix}^{\top}, \quad r = \rank(\mathbf{A}),
\end{equation}
where $
\begin{bmatrix} \widehat{\mathbf{U}}_{r} & \widehat{\mathbf{U}}_{c}
\end{bmatrix}$ and $
\begin{bmatrix} \widehat{\mathbf{V}}_{r} & \widehat{\mathbf{V}}_{c}
\end{bmatrix}$ are orthonormal, $\widehat{\mathbf{U}}_{r} \in \R^{m \times
r}$, $\widehat{\mathbf{V}}_{r} \in \R^{n \times r}$,
$\widehat{\mathbf{U}}_{c} \in \R^{m \times (m - r)}$,
$\widehat{\mathbf{V}}_{c} \in \R^{n \times (n - r)}$, and
$\mathbf{\Sigma}_{r} = \diag(\sigma_{1}, \ldots, \sigma_{r})$ with
$\sigma_{1} \ge \cdots \ge \sigma_{r} > 0$. The corresponding
singular vectors are
$\widehat{\mathbf{u}}_{i} = \widehat{\mathbf{U}}_{r} \mathbf{e}_{i}$ and
$\widehat{\mathbf{v}}_{i} = \widehat{\mathbf{V}}_{r} \mathbf{e}_{i}$.
The decomposition \eqref{eq:svd} yields the eigenvalue
decompositions
\begin{align}
    \mathbf{K}_{0}
    & =
    \begin{bmatrix}
        \mathbf{0} & \mathbf{A} \\
        \mathbf{A}^{\top} & \mathbf{0}
    \end{bmatrix}
    = \widehat{\mathbf{X}} \bm{\Lambda}_{0}
    \widehat{\mathbf{X}}^{\top}
    = \widehat{\mathbf{X}}
    \begin{bmatrix}
        \mathbf{\Sigma}_{r} & \mathbf{0} & \mathbf{0} \\
        \mathbf{0} & - \mathbf{\Sigma}_{r} & \mathbf{0} \\
        \mathbf{0} & \mathbf{0} & \mathbf{0}_{m + n - 2r}
    \end{bmatrix}
    \widehat{\mathbf{X}}^{\top}, \label{eq:k0-eigendecomposition} \\
    \shortintertext{and}
    \mathbf{K}^{2}
    & = \mathbf{I} + \mathbf{K}_{0}^{2}
    = \widehat{\mathbf{X}} \bm{\Lambda}
    \widehat{\mathbf{X}}^{\top}
    = \widehat{\mathbf{X}}
    \begin{bmatrix}
        \mathbf{I}_{r} + \mathbf{\Sigma}_{r}^{2} & \mathbf{0} & \mathbf{0} \\
        \mathbf{0} & \mathbf{I}_{r} + \mathbf{\Sigma}_{r}^{2} & \mathbf{0} \\
        \mathbf{0} & \mathbf{0} & \mathbf{I}_{m + n - 2r}
    \end{bmatrix}
    \widehat{\mathbf{X}}^{\top}, \label{eq:spectral-decompositions}
\end{align}
where
\[
    \widehat{\mathbf{X}} =
    \begin{bmatrix}
        \frac{1}{\sqrt{2}} \widehat{\mathbf{U}}_{r} &
        \frac{1}{\sqrt{2}} \widehat{\mathbf{U}}_{r} &
        \widehat{\mathbf{U}}_{c} & \mathbf{0} \\
        \frac{1}{\sqrt{2}} \widehat{\mathbf{V}}_{r} &
        - \frac{1}{\sqrt{2}} \widehat{\mathbf{V}}_{r} &
        \mathbf{0} & \widehat{\mathbf{V}}_{c}
    \end{bmatrix}
\]
is orthonormal.

\subsection{Properties of the block Krylov spaces}
\label{subsec:block-krylov}
For the case $\mathbf{M} = \mathbf{I}$ and $\mathbf{N} = \mathbf{I}$,
\Cref{alg:orth_tri} gives
\cite[Section 2.2]{montoison2021tricg}
\[
    \range(\mathbf{W}_{k})
    = \Kry_{k}(\mathbf{K}_{0}, \mathbf{B}_{0})
    = \spn \{ \mathbf{B}_{0}, \mathbf{K}_{0} \mathbf{B}_{0}, \ldots,
    \mathbf{K}_{0}^{k - 1} \mathbf{B}_{0} \}.
\]
Define
\[
    \mathbf{J} =
    \begin{bmatrix}
        \mathbf{I}_{m} & \mathbf{0} \\
        \mathbf{0} & - \mathbf{I}_{n}
    \end{bmatrix}.
\]
The block-diagonal structure of $\mathbf{W}_{k}$ in \eqref{eq:searchspace}
also gives $\range(\mathbf{J} \mathbf{W}_{k}) = \range(\mathbf{W}_{k})$.
With $\mathbf{J}_{2} = \diag(1, - 1)$, we have $\mathbf{J} \mathbf{B}_{0} =
\mathbf{B}_{0} \mathbf{J}_{2}$. Any $\mathbf{x} \in
\Kry_{k}(\mathbf{K}, \mathbf{B}_{0})$ can be represented as
\[
    \mathbf{x} = \sum_{\ell = 0}^{k - 1} \mathbf{K}^{\ell} \mathbf{B}_{0}
    \mathbf{a}_{\ell} =
    \sum_{s = 0}^{\lfloor(k - 1)/2 \rfloor}
    \mathbf{K}^{2s} \mathbf{B}_{0} \mathbf{a}_{2s}
    + \sum_{s = 0}^{\lfloor(k - 2)/2 \rfloor}
    \mathbf{K}^{2s + 1} \mathbf{B}_{0} \mathbf{a}_{2s + 1},
\]
where $\mathbf{a}_{\ell} \in \R^{2}$. Since $\mathbf{K}^{2} = \mathbf{I} +
\mathbf{K}_{0}^{2}$ and
$\mathbf{K} = \mathbf{K}_{0} + \mathbf{J}$,
\[
    \mathbf{K}^{2s} \mathbf{B}_{0}
    = (\mathbf{I} + \mathbf{K}_{0}^{2})^{s} \mathbf{B}_{0}, \quad
    \mathbf{K}^{2s + 1} \mathbf{B}_{0}
    = (\mathbf{I} + \mathbf{K}_{0}^{2})^{s}
    (\mathbf{K}_{0} \mathbf{B}_{0} + \mathbf{B}_{0} \mathbf{J}_{2}),
    \quad s \ge 0.
\]
Therefore, $\mathbf{K}^{\ell} \mathbf{B}_{0} \mathbf{a}_{\ell}$ is a linear
combination of the columns of
$\mathbf{B}_{0}, \mathbf{K}_{0} \mathbf{B}_{0}, \ldots, \mathbf{K}_{0}^{\ell}
\mathbf{B}_{0}$.
Hence $\mathbf{x} \in \Kry_{k}(\mathbf{K}_{0}, \mathbf{B}_{0})$.
Conversely, take any $\mathbf{x} \in \Kry_{k}(\mathbf{K}_{0}, \mathbf{B}_{0})$.
Since
\[
    \mathbf{K}_{0}^{2} = \mathbf{K}^{2} - \mathbf{I}, \quad
    \mathbf{K}_{0} \mathbf{B}_{0} = \mathbf{K} \mathbf{B}_{0} -
    \mathbf{B}_{0} \mathbf{J}_{2},
\]
the same argument gives $\mathbf{x} \in \Kry_{k}(\mathbf{K}, \mathbf{B}_{0})$.
The above discussion establishes the equivalence
\begin{equation}
    \label{eq:allspaces}
    \range(\mathbf{W}_{k}) = \Kry_{k}(\mathbf{K}_{0}, \mathbf{B}_{0}) =
    \Kry_{k}(\mathbf{K}, \mathbf{B}_{0}).
\end{equation}
Define
\begin{equation}
    \label{eq:startingblocks}
    \mathbf{D} =
    \begin{bmatrix}
        \mathbf{b} & \mathbf{A} \mathbf{c} & \mathbf{0} &
        \mathbf{0} \\
        \mathbf{0} & \mathbf{0} & - \mathbf{c} & \mathbf{A}^{\top}
        \mathbf{b}
    \end{bmatrix}.
\end{equation}
Then, we have $\range(\mathbf{D}) = \range\bigl(
    \begin{bmatrix} \mathbf{B}_{0} & \mathbf{K}_{0} \mathbf{B}_{0}
\end{bmatrix}\bigr)$.
For $k=2j$, by $\mathbf{K}^{2} = \mathbf{I} + \mathbf{K}_{0}^{2}$, we obtain
\begin{equation*}
    \begin{aligned}
        \range(\mathbf{W}_{2 j}) & = \Kry_{j}(\mathbf{K}_{0}^{2},
        \mathbf{B}_{0}) +
        \Kry_{j}(\mathbf{K}_{0}^{2}, \mathbf{K}_{0} \mathbf{B}_{0}) =
        \Kry_{j} \bigl(\mathbf{K}_{0}^{2},
            \begin{bmatrix}
                \mathbf{B}_{0} & \mathbf{K}_{0} \mathbf{B}_{0}
        \end{bmatrix} \bigr) \\
        & =
        \Kry_{j}(\mathbf{K}_{0}^{2}, \mathbf{D}) =
        \Kry_{j}(\mathbf{I} + \mathbf{K}_{0}^{2}, \mathbf{D}) \\
        & = \Kry_{j}(\mathbf{K}^{2}, \mathbf{D}).
    \end{aligned}
\end{equation*}
For $k = 2j+1$, we have
\begin{equation*}
    \range(\mathbf{W}_{2 j + 1}) = \Kry_{j}(\mathbf{K}^{2}, \mathbf{D}) +
    \range(\mathbf{K}^{2 j} \mathbf{B}_{0}).
\end{equation*}

\subsection{Relation between TriCG error and TriMR residual norms}
For an approximate solution $\mathbf{s}_{k}$, define the error at iteration
$k$ by
\[
    \bm{\varepsilon}_{k} := \mathbf{s}_{\star} - \mathbf{s}_{k}, \quad
    \mathbf{s}_{\star} = \mathbf{K}^{-1} \mathbf{f}
    =
    \begin{bmatrix}
        \mathbf{x}_{\star} \\
        \mathbf{y}_{\star}
    \end{bmatrix}.
\]
The superscripts $\cdot^\TCG$ and $\cdot^\TMR$ distinguish the two methods.
For either method, $\mathbf{r}_{k} := \mathbf{f} - \mathbf{K} \mathbf{s}_{k}
= \mathbf{K} \bm{\varepsilon}_{k}$. Since both methods start from zero,
they share the initial error and residual
\[
    \bm{\varepsilon}_{0} = \mathbf{s}_{\star}, \quad
    \mathbf{r}_{0} = \mathbf{f} = \mathbf{K} \bm{\varepsilon}_{0},
\]
for which the method superscripts are omitted.

The following theorem bounds the TriCG error norm by the TriMR
residual norm.

\begin{theorem} \label{thm:transfer}
    At the $k$th iteration, the TriCG error and TriMR residual norms satisfy
    \begin{equation}
        \label{eq:commoncomparison}
        \lVert \bm{\varepsilon}_{k}^{\TCG} \rVert
        \le \min_{\mathbf{v} \in \range(\mathbf{W}_{k})} \lVert
        \bm{\varepsilon}_{0} -
        \mathbf{v} \rVert_{\mathbf{K}^{2}}
        = \lVert \mathbf{r}_{k}^{\TMR} \rVert.
    \end{equation}
\end{theorem}

\begin{proof}
    From the block structures of $\mathbf{K}$ and $\mathbf{J}$, it follows that
    \begin{equation*}
        \mathbf{z}^{\top} \mathbf{J} \mathbf{K} \mathbf{z} =
        \lVert \mathbf{z} \rVert^{2}, \quad
        \lVert \mathbf{K} \mathbf{J} \mathbf{z} \rVert =
        \lVert \mathbf{z} \rVert_{\mathbf{K}^{2}}.
    \end{equation*}
    Since TriCG satisfies the Galerkin condition
    \eqref{eq:tricg_galerkin}, for any $\mathbf{v} \in
    \range(\mathbf{W}_{k})$, we have
    \[
        \mathbf{v}^{\top} \mathbf{r}_{k}^{\TCG} = \mathbf{v}^{\top}
        \mathbf{K} \bm{\varepsilon}_{k}^{\TCG} = 0.
    \]
    The relation $\range(\mathbf{J} \mathbf{W}_{k}) =
    \range(\mathbf{W}_{k})$ yields
    $\mathbf{v}^{\top} \mathbf{J} \mathbf{K} \bm{\varepsilon}_{k}^{\TCG} = 0$.
    Since $\mathbf{v} - \mathbf{s}_{k}^{\TCG} \in \range(\mathbf{W}_{k})$,
    \[
        (\mathbf{v} - \mathbf{s}_{k}^{\TCG})^{\top} \mathbf{JK}
        \bm{\varepsilon}_{k}^{\TCG} = 0.
    \]
    Hence,
    \[
        \begin{split}
            \lVert \bm{\varepsilon}_{k}^{\TCG} \rVert^{2}
            & = (\bm{\varepsilon}_{k}^{\TCG})^{\top} \mathbf{J} \mathbf{K}
            \bm{\varepsilon}_{k}^{\TCG} \\
            & = (\bm{\varepsilon}_{0} - \mathbf{v})^{\top} \mathbf{J}
            \mathbf{K} \bm{\varepsilon}_{k}^{\TCG}
            + (\mathbf{v} - \mathbf{s}_{k}^{\TCG})^{\top} \mathbf{J} \mathbf{K}
            \bm{\varepsilon}_{k}^{\TCG} \\
            & = (\mathbf{K} \mathbf{J} (\bm{\varepsilon}_{0} -
            \mathbf{v}))^{\top}
            \bm{\varepsilon}_{k}^{\TCG} \\
            & \le \lVert \bm{\varepsilon}_{0} - \mathbf{v}
            \rVert_{\mathbf{K}^{2}} \lVert
            \bm{\varepsilon}_{k}^{\TCG} \rVert.
        \end{split}
    \]
    Minimizing over $\mathbf{v}$ proves the first inequality. By the
    residual minimization property \eqref{eq:trimrsmall} of TriMR, we have
    \[
        \lVert \mathbf{r}_{k}^{\TMR} \rVert = \min_{\mathbf{v} \in
        \range(\mathbf{W}_{k})} \lVert \mathbf{f} - \mathbf{Kv} \rVert
        = \min_{\mathbf{v} \in \range(\mathbf{W}_{k})} \lVert
        \mathbf{K}(\bm{\varepsilon}_{0} - \mathbf{v}) \rVert =
        \min_{\mathbf{v} \in \range(\mathbf{W}_{k})} \lVert
        \bm{\varepsilon}_{0} - \mathbf{v} \rVert_{\mathbf{K}^{2}}. \qedhere
    \]
\end{proof}

By \Cref{thm:transfer}, any upper bound on
$\min_{\mathbf{v} \in \range(\mathbf{W}_{k})}
\lVert \bm{\varepsilon}_{0} - \mathbf{v} \rVert_{\mathbf{K}^{2}}$
also bounds both the TriCG error norm and the TriMR residual norm.
Our subsequent convergence analysis therefore focuses on bounding this quantity.

\section{Convergence bounds for even iterations}
\label{sec:even}
We derive convergence bounds at even iterations by estimating the
minimum in \eqref{eq:commoncomparison}. We first consider the full
spectral interval $[\lambda_{\min}(\mathbf{K}^{2}),
\lambda_{\max}(\mathbf{K}^{2})]$, then exploit the block
structures of $\mathbf{K}$ and
$\mathbf{D}$ to obtain bounds on smaller intervals.

\subsection{CG-type convergence bound and its refinement}
The following theorem gives convergence bounds analogous to the
classical CG bound \cite{trefethen1997numerical,saad2003iterative}.

\begin{theorem}
    \label{thm:full-spectrum}
    Consider TriCG and TriMR applied to \eqref{eq:identity-system} in
    exact arithmetic with a zero initial guess. Let
    $\mu = \lambda_{\min}(\mathbf{K}^{2})$ and
    $L = \lambda_{\max}(\mathbf{K}^{2})$.
    At every even iteration $k = 2j$, $j\ge1$, the TriCG error norm satisfies
    \begin{equation} \label{eq:cgunconditional}
        \lVert \bm{\varepsilon}_{2j}^{\TCG} \rVert \le
        2 \lVert \bm{\varepsilon}_{0} \rVert_{\mathbf{K}^{2}}
        \biggl(\frac{\sqrt{L/ \mu} - 1}{\sqrt{L/ \mu} + 1} \biggr)^{j},
    \end{equation}
    and the TriMR residual norm satisfies
    \begin{equation} \label{eq:mrunconditional}
        \lVert \mathbf{r}_{2j}^{\TMR} \rVert
        \le 2 \lVert \mathbf{r}_{0} \rVert
        \biggl(\frac{\sqrt{L/ \mu} - 1}{\sqrt{L/ \mu} + 1} \biggr)^{j}.
    \end{equation}
\end{theorem}

\begin{proof}
    Let $\mathbb{P}_{j}$ be the set of all polynomials of degree at
    most $j$. For any $p_{j} \in \mathbb{P}_{j}$ with $p_{j}(0) = 1$,
    we have $p_{j}(t) = 1 - t s_{j - 1}(t)$ for a polynomial
    $s_{j - 1}(t) \in \mathbb{P}_{j-1}$. The definition
    \eqref{eq:startingblocks} gives
    $\mathbf{K}^{2}\bm{\varepsilon}_{0}=\mathbf{K}\mathbf{r}_{0}
    =\mathbf{D}\one_{4}$. Thus the vector
    \begin{equation*} 
        \mathbf{v}
        = (\mathbf{I} - p_{j}(\mathbf{K}^{2})) \bm{\varepsilon}_{0}
        = s_{j - 1}(\mathbf{K}^{2}) \mathbf{K}^{2} \bm{\varepsilon}_{0}
        = s_{j - 1}(\mathbf{K}^{2}) \mathbf{D} \one_{4}
        \in \Kry_{j}(\mathbf{K}^{2}, \mathbf{D}) = \range(\mathbf{W}_{2j}).
    \end{equation*}
    Hence
    \begin{equation}
        \label{eq:polynomial-approximation-bound}
        \min_{\mathbf{w} \in \range(\mathbf{W}_{2j})}
        \lVert \bm{\varepsilon}_{0} - \mathbf{w} \rVert_{\mathbf{K}^{2}}
        \le \lVert \bm{\varepsilon}_{0} - \mathbf{v} \rVert_{\mathbf{K}^{2}}
        = \lVert p_{j}(\mathbf{K}^{2}) \bm{\varepsilon}_{0}
        \rVert_{\mathbf{K}^{2}}.
    \end{equation}
    To bound the right-hand side, we use the decomposition
    \eqref{eq:spectral-decompositions} to obtain
    \begin{align*}
        \lVert p_{j}(\mathbf{K}^{2}) \bm{\varepsilon}_{0}
        \rVert_{\mathbf{K}^{2}}^{2}
        & = \bigl\lVert \bm{\Lambda}^{1/2}
        p_{j}(\bm{\Lambda})
        \widehat{\mathbf{X}}^{\top} \bm{\varepsilon}_{0} \bigr\rVert^{2}
        \notag \\
        & \le \max_{\lambda \in [\mu, L]} \lvert p_{j}(\lambda) \rvert^{2}
        \bigl\lVert \bm{\Lambda}^{1/2}
        \widehat{\mathbf{X}}^{\top} \bm{\varepsilon}_{0} \bigr\rVert^{2}
        \notag \\
        & = \max_{\lambda \in [\mu, L]} \lvert p_{j}(\lambda) \rvert^{2}
        \lVert \bm{\varepsilon}_{0} \rVert_{\mathbf{K}^{2}}^{2}.
    \end{align*}
    Let $T_{j}$ be the Chebyshev polynomial of the first kind. In particular, we
    consider the polynomial
    \begin{equation*} 
        p_{j}(t) =
        T_{j} \Bigl(\frac{L + \mu - 2t}{L - \mu} \Bigr)
        \Big/
        T_{j} \Bigl(\frac{L + \mu}{L - \mu} \Bigr).
    \end{equation*}
    This polynomial has degree at most $j$ and satisfies $p_{j}(0) = 1$.
    For $t \in [\mu, L]$, the argument $(L + \mu - 2t)/(L - \mu)$ belongs
    to $[-1, 1]$, where $\lvert T_{j} \rvert \le 1$. 
    For $z > 1$, the Chebyshev polynomial satisfies
    \begin{equation}\label{eq:chebyshev-explicit}
        T_{j}(z) = \frac{1}{2}
        \bigl[
            \bigl(z + \sqrt{z^{2} - 1} \bigr)^{j}
            + \bigl(z - \sqrt{z^{2} - 1} \bigr)^{j}
        \bigr].
    \end{equation}
    Therefore,
    \[
        T_{j} \Bigl(\frac{L + \mu}{L - \mu} \Bigr)
        = \frac{1}{2}
        \biggl[
            \biggl(\frac{\sqrt{L/ \mu} + 1}{\sqrt{L/ \mu} - 1} \biggr)^{j}
            + \biggl(\frac{\sqrt{L/ \mu} + 1}{\sqrt{L/ \mu} - 1} \biggr)^{-j}
        \biggr],
    \]
    and
    \begin{equation}
        \label{eq:chebbound}
        \max_{t \in [\mu, L]} \lvert p_{j}(t) \rvert
        \le T_{j} \Bigl(\frac{L + \mu}{L - \mu} \Bigr)^{-1} \le 2
        \biggl(\frac{\sqrt{L/ \mu} - 1}{\sqrt{L/ \mu} + 1} \biggr)^{j} .
    \end{equation}
    From \Cref{thm:transfer} and
    \eqref{eq:polynomial-approximation-bound}--\eqref{eq:chebbound},
    \eqref{eq:cgunconditional} and \eqref{eq:mrunconditional} hold.
\end{proof}

The convergence bounds \eqref{eq:cgunconditional} and
\eqref{eq:mrunconditional} are analogous to the classical CG convergence
bound. When $\Null(\mathbf{K}_{0})$ is nonempty, $\mu = 1$.
We next show that the convergence factor can be improved
by replacing $L/ \mu$ with $L/(1 + \sigma_{r}^{2})$.

For the symmetric matrix $\mathbf{K}_{0}$ defined in
\eqref{eq:k0-eigendecomposition}, let $\mathbf{P}_{\mathbf{K}_{0}}$ and
$\mathbf{P}_{\mathbf{K}_{0}}^{\perp}$ be orthogonal projectors onto
$\range(\mathbf{K}_{0})$ and $\Null(\mathbf{K}_{0})$, respectively.
From the eigenvalue decomposition \eqref{eq:k0-eigendecomposition} of
$\mathbf{K}_{0}$, it follows that
\begin{equation} \label{eq:projectors}
    \mathbf{P}_{\mathbf{K}_{0}}
    =
    \begin{bmatrix}
        \widehat{\mathbf{U}}_{r} \widehat{\mathbf{U}}_{r}^{\top} & \mathbf{0} \\
        \mathbf{0} & \widehat{\mathbf{V}}_{r} \widehat{\mathbf{V}}_{r}^{\top}
    \end{bmatrix}, \quad
    \mathbf{P}_{\mathbf{K}_{0}}^{\perp} = \mathbf{I} -
    \mathbf{P}_{\mathbf{K}_{0}}
    =
    \begin{bmatrix}
        \widehat{\mathbf{U}}_{c} \widehat{\mathbf{U}}_{c}^{\top} & \mathbf{0} \\
        \mathbf{0} & \widehat{\mathbf{V}}_{c} \widehat{\mathbf{V}}_{c}^{\top}
    \end{bmatrix}.
\end{equation}
The following theorem can improve the convergence factor presented in
\Cref{thm:full-spectrum}.

\begin{theorem} \label{thm:null}
    Consider TriCG and TriMR applied to \eqref{eq:identity-system} in
    exact arithmetic with a zero initial guess.
    Suppose $0<\sigma_{r}<\sigma_{1}$.
    At every even iteration $k = 2j$, $j\ge1$, the TriCG error norm satisfies
    \begin{equation} \label{eq:bound2}
        \lVert \bm{\varepsilon}_{2j}^{\TCG} \rVert
        \le C \biggl(
            \frac{\sqrt{(1 + \sigma_{1}^{2})/(1 + \sigma_{r}^{2})} - 1}
            {\sqrt{(1 + \sigma_{1}^{2})/(1 + \sigma_{r}^{2})} + 1}
        \biggr)^{j},
    \end{equation}
    where
    \[
        C=2\Bigl(
            \|\mathbf{P}_{\mathbf{K}_{0}}\bm{\varepsilon}_{0}\|_{\mathbf{K}^{2}}
            +\frac{\sigma_{1}+\sigma_{r}}{\sigma_{1}-\sigma_{r}}
            \|\mathbf{P}_{\mathbf{K}_{0}}\mathbf{J}\mathbf{r}_{0}\|_{\mathbf{K}^{2}}
        \Bigr).
    \]
    The same bound holds for $\lVert \mathbf{r}_{2j}^{\TMR} \rVert$.
\end{theorem}
\begin{proof}
    Let $\tau_{1}=1+\sigma_{r}^{2}$ and $L=1+\sigma_{1}^{2}$.
    For $\ell\ge0$, define
    \[
        \eta(t)=\frac{L+\tau_{1}-2t}{L-\tau_{1}},\quad
        \gamma=\eta(1)
        =\frac{\sigma_{1}^{2}+\sigma_{r}^{2}}{\sigma_{1}^{2}-\sigma_{r}^{2}}>1,
        \quad
        \phi_{\ell}(t)=\frac{T_{\ell}(\eta(t))}{T_{\ell}(\gamma)}.
    \]
    By construction, $\phi_{\ell}(1)=1$. Moreover,
    $\eta(0)>\gamma>1$ implies $\phi_{\ell}(0)>0$.

    From the block structures of $\mathbf{P}_{\mathbf{K}_{0}}$,
    $\mathbf{P}_{\mathbf{K}_{0}}^{\perp}$, and $\mathbf{J}$, we have
    \begin{equation} \label{eq:JP}
        \mathbf{P}_{\mathbf{K}_{0}} \mathbf{J} = \mathbf{J}
        \mathbf{P}_{\mathbf{K}_{0}}, \quad
        \mathbf{P}_{\mathbf{K}_{0}}^{\perp}
        \mathbf{J} = \mathbf{J} \mathbf{P}_{\mathbf{K}_{0}}^{\perp}.
    \end{equation}
    Since $\mathbf{P}_{\mathbf{K}_{0}} \mathbf{K}_{0} = \mathbf{K}_{0}
    \mathbf{P}_{\mathbf{K}_{0}} = \mathbf{K}_{0}$ and
    $\mathbf{P}_{\mathbf{K}_{0}}^{\perp} \mathbf{K}_{0} = \mathbf{0}$, we have
    \begin{equation} \label{eq:PK}
        \mathbf{P}_{\mathbf{K}_{0}} \mathbf{K}
        =\mathbf{P}_{\mathbf{K}_{0}} (\mathbf{J} + \mathbf{K}_{0})
        =(\mathbf{J} + \mathbf{K}_{0}) \mathbf{P}_{\mathbf{K}_{0}}
        =\mathbf{K} \mathbf{P}_{\mathbf{K}_{0}},\quad
        \mathbf{P}_{\mathbf{K}_{0}}^{\perp} \mathbf{K}
        =\mathbf{P}_{\mathbf{K}_{0}}^{\perp} (\mathbf{J} + \mathbf{K}_{0})
        =\mathbf{J} \mathbf{P}_{\mathbf{K}_{0}}^{\perp}.
    \end{equation}
    Then, by \eqref{eq:JP} and \eqref{eq:PK}, we obtain
    \begin{equation}
        \mathbf{P}_{\mathbf{K}_{0}}^{\perp} \mathbf{J}\mathbf{r}_{0}
        =\mathbf{P}_{\mathbf{K}_{0}}^{\perp} \mathbf{J}\mathbf{K}
        \bm{\varepsilon}_{0}
        = \mathbf{J}^{2} \mathbf{P}_{\mathbf{K}_{0}}^{\perp}
        \bm{\varepsilon}_{0}
        =\mathbf{P}_{\mathbf{K}_{0}}^{\perp} \bm{\varepsilon}_{0}.
        \label{eq:matching}
    \end{equation}
    Construct the vector
    \[
        \mathbf{v} = \biggl(\mathbf{I}-
        \frac{\phi_{j}(\mathbf{K}^{2})}{\phi_{j}(0)}\biggr) \bm{\varepsilon}_{0}
        +\frac{\phi_{j-1}(\mathbf{K}^{2})}{\phi_{j}(0)}\mathbf{J}\mathbf{r}_{0}.
    \]
    By the same argument as in the proof of \Cref{thm:full-spectrum},
    the first term on the right-hand side belongs to
    $\Kry_{j}(\mathbf{K}^{2}, \mathbf{D}) = \range(\mathbf{W}_{2j})$.
    From the definition \eqref{eq:startingblocks} of $\mathbf{D}$, we
    have $\mathbf{J}\mathbf{r}_{0} = \mathbf{D}(\mathbf{e}_{1} +
    \mathbf{e}_{3}) \in
    \range(\mathbf{D})$. Hence,
    \[
        \mathbf{v} \in \Kry_{j}(\mathbf{K}^{2}, \mathbf{D}) =
        \range(\mathbf{W}_{2j}).
    \]
    It follows from \eqref{eq:PK} that
    \begin{equation} \label{eq:PK2}
        \mathbf{P}_{\mathbf{K}_{0}} \mathbf{K}^{2} = \mathbf{K}
        \mathbf{P}_{\mathbf{K}_{0}} \mathbf{K} = \mathbf{K}^{2}
        \mathbf{P}_{\mathbf{K}_{0}}, \quad
        \mathbf{P}_{\mathbf{K}_{0}}^{\perp} \mathbf{K}^{2} = \mathbf{J}
        \mathbf{P}_{\mathbf{K}_{0}}^{\perp} \mathbf{K} = \mathbf{J}^{2}
        \mathbf{P}_{\mathbf{K}_{0}}^{\perp} =
        \mathbf{P}_{\mathbf{K}_{0}}^{\perp}.
    \end{equation}
    Then, by $\phi_{\ell}(1) = 1$ and
    \eqref{eq:matching}--\eqref{eq:PK2}, we obtain
    \[
        \mathbf{P}_{\mathbf{K}_{0}}^{\perp} (\bm{\varepsilon}_{0} - \mathbf{v})
        = \frac{\phi_{j}(1)}{\phi_{j}(0)}
        \mathbf{P}_{\mathbf{K}_{0}}^{\perp} \bm{\varepsilon}_{0}
        - \frac{\phi_{j-1}(1)}{\phi_{j}(0)}
        \mathbf{P}_{\mathbf{K}_{0}}^{\perp} \mathbf{J}\mathbf{r}_{0}
        = \mathbf{0},
    \]
    and
    \[
        \mathbf{P}_{\mathbf{K}_{0}} (\bm{\varepsilon}_{0} - \mathbf{v})
        = \frac{\phi_{j}(\mathbf{K}^{2})}{\phi_{j}(0)}
        \mathbf{P}_{\mathbf{K}_{0}} \bm{\varepsilon}_{0}
        - \frac{\phi_{j-1}(\mathbf{K}^{2})}{\phi_{j}(0)}
        \mathbf{P}_{\mathbf{K}_{0}} \mathbf{J}\mathbf{r}_{0}.
    \]
    Hence,
    \begin{align}
        \min_{\mathbf{w} \in \range(\mathbf{W}_{2j})} \lVert
        \bm{\varepsilon}_0 - \mathbf{w} \rVert_{\mathbf{K}^2}^2
        &\le \lVert \bm{\varepsilon}_0 - \mathbf{v}
        \rVert_{\mathbf{K}^2}^2 \notag \\
        &= \lVert \mathbf{P}_{\mathbf{K}_0} (\bm{\varepsilon}_0 -
        \mathbf{v}) + \mathbf{P}_{\mathbf{K}_0}^\perp
        (\bm{\varepsilon}_0 - \mathbf{v}) \rVert_{\mathbf{K}^2}^2 \notag \\
        &= \lVert \mathbf{P}_{\mathbf{K}_0} (\bm{\varepsilon}_0 -
        \mathbf{v}) \rVert_{\mathbf{K}^2}^2. \label{eq:Pkerr}
    \end{align}
    It remains to bound
    $\|\mathbf{P}_{\mathbf{K}_{0}}(\bm{\varepsilon}_{0}-\mathbf{v})\|_{\mathbf{K}^{2}}$.

    Since $\eta(t)\in[-1,1]$ for $t \in [\tau_{1},L]$,
    $|\phi_{\ell}(t)|\le T_{\ell}(\gamma)^{-1}$.
    By the eigenvalue decomposition
    \eqref{eq:spectral-decompositions} and \eqref{eq:projectors},
    we have
    \[
        \|\phi_{\ell}(\mathbf{K}^{2})\mathbf{P}_{\mathbf{K}_{0}}\|
        \le\max_{t\in[\tau_{1},L]}|\phi_{\ell}(t)|
        \le\frac{1}{T_{\ell}(\gamma)}.
    \]
    For $\mathbf{x}\in\range(\mathbf{K}_{0})$, we have
    $\mathbf{P}_{\mathbf{K}_{0}}\mathbf{x}=\mathbf{x}$.
    Using \eqref{eq:PK} and $\mathbf{K} \phi_{\ell}(\mathbf{K}^2) =
    \phi_{\ell}(\mathbf{K}^2) \mathbf{K}$, we obtain
    \begin{align}
        \|\phi_{\ell}(\mathbf{K}^{2})\mathbf{x}\|_{\mathbf{K}^{2}}
        &=\|\mathbf{K}\phi_{\ell}(\mathbf{K}^{2})
        \mathbf{P}_{\mathbf{K}_{0}} \mathbf{x}\|
        =\|\phi_{\ell}(\mathbf{K}^{2})\mathbf{P}_{\mathbf{K}_{0}}
        \mathbf{K}\mathbf{x}\|\notag \\
        &\le\|\phi_{\ell}(\mathbf{K}^{2})\mathbf{P}_{\mathbf{K}_{0}}\|
        \|\mathbf{K}\mathbf{x}\| \notag \\
        &\le \max_{t\in [\tau_{1}, L]} \lvert \phi_{\ell}(t) \rvert \lVert
        \mathbf{x} \rVert_{\mathbf{K}^{2}}
        \le\frac{\|\mathbf{x}\|_{\mathbf{K}^{2}}}{T_{\ell}(\gamma)}.
        \label{eq:weighted-polynomial-estimate}
    \end{align}
    Both $\mathbf{P}_{\mathbf{K}_{0}}\bm{\varepsilon}_{0}$ and
    $\mathbf{P}_{\mathbf{K}_{0}}\mathbf{J}\mathbf{r}_{0}$ belong to
    $\range(\mathbf{K}_{0})$. Applying the triangle inequality and
    \eqref{eq:weighted-polynomial-estimate} with $\ell=j$ and $\ell=j-1$,
    together with $\phi_{j}(0)=T_{j}(\eta(0))/T_{j}(\gamma)$, yields
    \begin{align*}
        \|\mathbf{P}_{\mathbf{K}_{0}}(\bm{\varepsilon}_{0}-\mathbf{v})\|_{\mathbf{K}^{2}}
        &=\biggl\|
        \frac{\phi_{j}(\mathbf{K}^{2})}{\phi_{j}(0)}
        \mathbf{P}_{\mathbf{K}_{0}}\bm{\varepsilon}_{0}
        -\frac{\phi_{j-1}(\mathbf{K}^{2})}{\phi_{j}(0)}
        \mathbf{P}_{\mathbf{K}_{0}}\mathbf{J}\mathbf{r}_{0}
        \biggr\|_{\mathbf{K}^{2}}\\
        &\le\frac{1}{\phi_{j}(0)} \bigl(
            \|\phi_{j}(\mathbf{K}^{2})\mathbf{P}_{\mathbf{K}_{0}}
            \bm{\varepsilon}_{0}\|_{\mathbf{K}^{2}}
            +\|\phi_{j-1}(\mathbf{K}^{2})\mathbf{P}_{\mathbf{K}_{0}}
        \mathbf{J}\mathbf{r}_{0}\|_{\mathbf{K}^{2}}\bigr)\\
        &\le\frac{T_{j}(\gamma)}{T_{j}(\eta(0))}\biggl(
            \frac{\|\mathbf{P}_{\mathbf{K}_{0}}\bm{\varepsilon}_{0}\|_{\mathbf{K}^{2}}}
            {T_{j}(\gamma)}
            +\frac{\|\mathbf{P}_{\mathbf{K}_{0}}\mathbf{J}\mathbf{r}_{0}\|_{\mathbf{K}^{2}}}
        {T_{j-1}(\gamma)}\biggr)\\
        &=\frac{1}{T_{j}(\eta(0))}\Bigl(
            \|\mathbf{P}_{\mathbf{K}_{0}}\bm{\varepsilon}_{0}\|_{\mathbf{K}^{2}}
            +\frac{T_{j}(\gamma)}{T_{j-1}(\gamma)}
            \|\mathbf{P}_{\mathbf{K}_{0}}\mathbf{J}\mathbf{r}_{0}\|_{\mathbf{K}^{2}}
        \Bigr).
    \end{align*}
    To bound the coefficient independently of $j$, let
    \[
        \rho_{j} = \frac{T_j(\gamma)}{T_{j-1}(\gamma)}, \quad
        \rho=\gamma+\sqrt{\gamma^{2}-1}
        =\frac{\sigma_{1}+\sigma_{r}}{\sigma_{1}-\sigma_{r}} > 1.
    \]
    Since $\gamma-\sqrt{\gamma^{2}-1}=\rho^{-1}$,
    \eqref{eq:chebyshev-explicit} gives
    $T_{\ell}(\gamma)=(\rho^{\ell}+\rho^{-\ell})/2$. Hence,
    \begin{equation}\label{eq:chebyshev-ratio}
        \rho_{j}=\frac{T_{j}(\gamma)}{T_{j-1}(\gamma)}
        =\frac{\rho^{j}+\rho^{-j}}{\rho^{j-1}+\rho^{-j+1}}
        =\rho\frac{1+\rho^{-2j}}{1+ \rho^2 \rho^{-2j}}
        \le\rho.
    \end{equation}
    The last inequality follows from $\rho>1$.
    Using \eqref{eq:chebyshev-ratio} and the Chebyshev estimate obtained
    from \eqref{eq:chebbound}, we obtain
    \begin{align}
        \|\mathbf{P}_{\mathbf{K}_{0}}(\bm{\varepsilon}_{0}-\mathbf{v})\|_{\mathbf{K}^{2}}
        &\le T_{j}(\eta(0))^{-1} \bigl(
            \|\mathbf{P}_{\mathbf{K}_{0}}\bm{\varepsilon}_{0}\|_{\mathbf{K}^{2}}
            +\rho_{j}\|\mathbf{P}_{\mathbf{K}_{0}}\mathbf{J}\mathbf{r}_{0}\|_{\mathbf{K}^{2}}
        \bigr) \notag \\
        &=
        \bigl(\|\mathbf{P}_{\mathbf{K}_{0}}\bm{\varepsilon}_{0}\|_{\mathbf{K}^{2}}
        +\rho_{j}\|\mathbf{P}_{\mathbf{K}_{0}}\mathbf{J}\mathbf{r}_{0}\|_{\mathbf{K}^{2}}\bigr)
        T_j\Bigl(\frac{L+\tau_1}{L-\tau_1}\Bigr)^{-1} \notag \\
        &\le
        2\bigl(\|\mathbf{P}_{\mathbf{K}_{0}}\bm{\varepsilon}_{0}\|_{\mathbf{K}^{2}}
            +\rho
        \|\mathbf{P}_{\mathbf{K}_{0}}\mathbf{J}\mathbf{r}_{0}\|_{\mathbf{K}^{2}}\bigr)
        \biggl(
        \frac{\sqrt{L/\tau_{1}}-1}{\sqrt{L/\tau_{1}}+1}\biggr)^{j}.
        \label{eq:PKerr2}
    \end{align}
    Combining \eqref{eq:Pkerr} and \eqref{eq:PKerr2} with
    \Cref{thm:transfer} proves the stated
    bounds for the TriCG error and the TriMR residual norms.
\end{proof}

\subsection{Convergence bounds on a reduced spectral interval}
\label{subsec:blockcg}
In this subsection, we assume $r\ge3$ and derive sufficient conditions
for replacing $L/(1 + \sigma_{r}^{2})$ in the convergence factor of
\Cref{thm:null} with $L/(1 + \sigma_{r-1}^{2})$.
The resulting factor is no larger and is strictly smaller when
$\sigma_{r-1}>\sigma_{r}$. This does not necessarily yield a smaller
bound at every iteration, since the coefficients also differ.
Define
\[
    \mathbf{b}_{0} = \widehat{\mathbf{U}}_{c} \widehat{\mathbf{U}}_{c}^{\top}
    \mathbf{b}, \quad
    \mathbf{c}_{0} = \widehat{\mathbf{V}}_{c} \widehat{\mathbf{V}}_{c}^{\top}
    \mathbf{c}.
\]
Let the columns of $\mathbf{Z}_{0} \in \R^{(m+n) \times d_{0}}
\ (d_{0} \le 2)$ form an orthonormal basis for
\begin{equation*}
    \range\biggl(
        \begin{bmatrix}
            \mathbf{b}_{0} & \mathbf{0} \\
            \mathbf{0} & \mathbf{c}_{0}
    \end{bmatrix}\biggr).
\end{equation*}
By \eqref{eq:projectors}, \eqref{eq:matching}, and
$\mathbf{r}_{0}=\mathbf{f}$, we obtain
\begin{equation}
    \label{eq:nullerror}
    \mathbf{P}_{\mathbf{K}_{0}}^{\perp} \bm{\varepsilon}_{0}
    = \mathbf{P}_{\mathbf{K}_{0}}^{\perp} \mathbf{J}\mathbf{f}
    =
    \begin{bmatrix}
        \mathbf{b}_{0} \\ -\mathbf{c}_{0}
    \end{bmatrix}
    \in \range(\mathbf{Z}_{0}).
\end{equation}
Moreover, by \eqref{eq:startingblocks} and \eqref{eq:projectors},
\begin{align*}
    \mathbf{P}_{\mathbf{K}_{0}}^{\perp} \mathbf{D}
    & =
    \begin{bmatrix}
        \widehat{\mathbf{U}}_{c} \widehat{\mathbf{U}}_{c}^{\top} & \mathbf{0} \\
        \mathbf{0} & \widehat{\mathbf{V}}_{c} \widehat{\mathbf{V}}_{c}^{\top}
    \end{bmatrix}
    \begin{bmatrix}
        \mathbf{b} & \mathbf{A} \mathbf{c} & \mathbf{0} & \mathbf{0} \\
        \mathbf{0} & \mathbf{0} & - \mathbf{c} & \mathbf{A}^{\top} \mathbf{b}
    \end{bmatrix} \notag \\
    & =
    \begin{bmatrix}
        \mathbf{b}_{0} &
        \widehat{\mathbf{U}}_{c} \widehat{\mathbf{U}}_{c}^{\top}
        \mathbf{A} \mathbf{c}
        & \mathbf{0} & \mathbf{0} \\
        \mathbf{0} & \mathbf{0} & - \mathbf{c}_{0} &
        \widehat{\mathbf{V}}_{c} \widehat{\mathbf{V}}_{c}^{\top}
        \mathbf{A}^{\top} \mathbf{b}
    \end{bmatrix} \notag \\
    & =
    \begin{bmatrix}
        \mathbf{b}_{0} & \mathbf{0} & \mathbf{0} & \mathbf{0} \\
        \mathbf{0} & \mathbf{0} & - \mathbf{c}_{0} & \mathbf{0}
    \end{bmatrix},
\end{align*}
where the last equality follows from
$\widehat{\mathbf{U}}_{c}^{\top} \mathbf{A} = \mathbf{0}$ and
$\widehat{\mathbf{V}}_{c}^{\top} \mathbf{A}^{\top} = \mathbf{0}$.
Hence,
\begin{equation} \label{eq:null-starting-block}
    \range(\mathbf{P}_{\mathbf{K}_0}^\perp \mathbf{D}) \subseteq
    \range(\mathbf{Z}_0).
\end{equation}

Define
\begin{equation}
    \label{eq:Xspecial}
    \mathbf{X}_{1} =
    \begin{bmatrix}
        \mathbf{Z}_{0} &
        \begin{matrix}
            \widehat{\mathbf{u}}_{r} & \mathbf{0} \\
            \mathbf{0} & \widehat{\mathbf{v}}_{r}
        \end{matrix}
    \end{bmatrix}, \quad
    \bm{\Lambda}_{1} =
    \begin{bmatrix}
        \mathbf{I}_{d_{0}} & \mathbf{0} & \mathbf{0} \\
        \mathbf{0} & 1 + \sigma_{r}^{2} & 0 \\
        \mathbf{0} & 0 & 1 + \sigma_{r}^{2}
    \end{bmatrix} \in \R^{d \times d}, \quad d = d_{0} + 2 \le 4,
\end{equation}
and
\[
    \mathbf{X}_{2} =
    \begin{bmatrix}
        \widehat{\mathbf{u}}_{1} & \mathbf{0} & \cdots &
        \widehat{\mathbf{u}}_{r-1} & \mathbf{0} \\
        \mathbf{0} & \widehat{\mathbf{v}}_{1} & \cdots & \mathbf{0} &
        \widehat{\mathbf{v}}_{r-1}
    \end{bmatrix}, \quad
    \bm{\Lambda}_{2} = \diag(1 + \sigma_{1}^{2}, 1 + \sigma_{1}^{2}, \ldots,
    1 + \sigma_{r-1}^{2}, 1 + \sigma_{r-1}^{2}).
\]
The columns of $
\begin{bmatrix}\mathbf{X}_{1} & \mathbf{X}_{2}
\end{bmatrix}$
are orthonormal and satisfy
\begin{equation}\label{eq:spectral-block-relations}
    \mathbf{K}^{2} \mathbf{X}_{1} = \mathbf{X}_{1}
    \bm{\Lambda}_{1}, \quad
    \mathbf{K}^{2} \mathbf{X}_{2} = \mathbf{X}_{2} \bm{\Lambda}_{2}.
\end{equation}
By \eqref{eq:projectors} and \eqref{eq:nullerror},
\begin{equation}\label{eq:initial-error-subspace}
    \bm{\varepsilon}_{0}
    = \mathbf{P}_{\mathbf{K}_{0}} \bm{\varepsilon}_{0}
    + \mathbf{P}_{\mathbf{K}_{0}}^{\perp} \bm{\varepsilon}_{0}
    \in \range \Biggl(
        \begin{bmatrix}
            \widehat{\mathbf{U}}_{r} & \mathbf{0} \\
            \mathbf{0} & \widehat{\mathbf{V}}_{r}
    \end{bmatrix} \Biggr)
    + \range(\mathbf{Z}_{0})
    = \range \bigl(
        \begin{bmatrix}
            \mathbf{X}_{1} & \mathbf{X}_{2}
    \end{bmatrix} \bigr).
\end{equation}
Similarly, \eqref{eq:null-starting-block} gives
\begin{equation}\label{eq:starting-block-subspace}
    \range(\mathbf{D}) \subseteq \range \bigl(
        \begin{bmatrix}
            \mathbf{X}_{1} & \mathbf{X}_{2}
    \end{bmatrix} \bigr).
\end{equation}
Since, for $0\le\ell\le j-1$,
\begin{align*}
    (\mathbf{K}^{2})^{\ell} \mathbf{D}
    & =  (\mathbf{K}^{2})^{\ell} (\mathbf{X}_{1}
        \mathbf{X}_{1}^{\top} \mathbf{D} + \mathbf{X}_{2}
        \mathbf{X}_{2}^{\top}
    \mathbf{D}) \\
    & = \mathbf{X}_{1} (\bm{\Lambda}_{1}^{\ell} \mathbf{X}_{1}^{\top}
    \mathbf{D}) + \mathbf{X}_{2} (\bm{\Lambda}_{2}^{\ell} \mathbf{X}_{2}^{\top}
    \mathbf{D}),
\end{align*}
we have
\begin{equation}\label{eq:krylov-spectral-subspace}
    \range(\mathbf{W}_{2j}) =
    \Kry_{j}(\mathbf{K}^{2}, \mathbf{D}) \subseteq \range \bigl(
        \begin{bmatrix}
            \mathbf{X}_{1} & \mathbf{X}_{2}
    \end{bmatrix} \bigr).
\end{equation}
\Needspace{6\baselineskip}
The following theorem can further improve the convergence factor in
\Cref{thm:null}.

\begin{theorem}\label{thm:cgconditional}
    For the construction in \eqref{eq:Xspecial}, assume that
    $0 < \sigma_{r-1} < \sigma_{1}$ and that $\mathbf{R}_{1} =
    \mathbf{X}_{1}^{\top} \mathbf{D}$ has full row rank. At every even
    iteration $k=2j$, $j\ge1$, the TriCG error norm satisfies
    \begin{equation}
        \label{eq:conditionalCG}
        \lVert \bm{\varepsilon}_{2j}^{\TCG} \rVert
        \le B \left(
            \frac{\sqrt{(1 + \sigma_{1}^{2})/(1 + \sigma_{r-1}^{2})} - 1}
            {\sqrt{(1 + \sigma_{1}^{2})/(1 + \sigma_{r-1}^{2})} + 1}
        \right)^{j},
    \end{equation}
    where
    \[
        B = 2\Bigl(
            \|\bm{\Lambda}_{2}^{1/2}\mathbf{X}_{2}^\top \bm{\varepsilon}_0\|
            +\frac{\sigma_{1}+\sigma_{r-1}}{\sigma_{1}-\sigma_{r-1}}
            \|\bm{\Lambda}_{2}^{1/2}(\mathbf{X}_{2}^{\top} \mathbf{D}
            \mathbf{R}_{1}^{\dagger})\|
        \,\|\mathbf{X}_{1}^\top \bm{\varepsilon}_0\|\Bigr).
    \]
    The same bound holds for $\lVert \mathbf{r}_{2j}^{\TMR} \rVert$.
\end{theorem}

\begin{proof}
    Let $\tau_{2} = 1 + \sigma_{r-1}^{2}$ and $L = 1 + \sigma_{1}^{2}$.
    From \eqref{eq:initial-error-subspace}, the initial error can be
    represented by
    \[
        \bm{\varepsilon}_{0} = \mathbf{X}_{1} \bm{\xi}_{1} +
        \mathbf{X}_{2} \bm{\xi}_{2}, \quad
        \bm{\xi}_{1} = \mathbf{X}_{1}^{\top} \bm{\varepsilon}_{0}, \quad
        \bm{\xi}_{2} = \mathbf{X}_{2}^{\top} \bm{\varepsilon}_{0}.
    \]
    Let $\mathbf{R}_{1}^{\dagger} \in \R^{4 \times d}$ be the
    Moore--Penrose pseudoinverse of $\mathbf{R}_{1} =
    \mathbf{X}_{1}^{\top} \mathbf{D}$.
    The full row rank condition implies
    $\mathbf{R}_{1} \mathbf{R}_{1}^{\dagger} = \mathbf{I}_{d}$. Define
    $\mathbf{F} = \mathbf{D} \mathbf{R}_{1}^{\dagger}$. From
    \eqref{eq:starting-block-subspace}, we have
    \begin{align*}
        \mathbf{F} = \mathbf{D} \mathbf{R}_{1}^{\dagger}
        = \bigl(\mathbf{X}_{1} (\mathbf{X}_{1}^{\top} \mathbf{D})
        + \mathbf{X}_{2} (\mathbf{X}_{2}^{\top} \mathbf{D}) \bigr)
        \mathbf{R}_{1}^{\dagger}
        = \mathbf{X}_{1} + \mathbf{X}_{2} \bm{\Gamma}, \quad
        \bm{\Gamma} = \mathbf{X}_{2}^{\top} \mathbf{D}
        \mathbf{R}_{1}^{\dagger}.
    \end{align*}
    Define
    \[
        \eta(t) = \frac{L + \tau_{2} - 2t}{L - \tau_{2}}, \quad
        \psi_{\ell}(t) = T_{\ell}(\eta(t)), \quad \ell\ge0.
    \]
    For each diagonal entry $\lambda_{i}(\bm{\Lambda}_{1})
    =\mathbf{e}_{i}^{\top}\bm{\Lambda}_{1}\mathbf{e}_{i}
    \le\tau_{2}$, $1\le i\le d$, we have
    \[
        \eta(\lambda_i(\bm{\Lambda}_{1}))
        = 1 + \frac{2(\tau_{2} - \lambda_i(\bm{\Lambda}_{1}))}{L -
        \tau_{2}} \ge 1.
    \]
    Formula~\eqref{eq:chebyshev-explicit} also holds at $x=1$, where
    $T_{\ell}(1)=1$. Thus $T_{\ell}(x)\ge1$ for $x\ge1$ and $\ell\ge0$.
    Hence, the diagonal matrix
    \[
        \psi_{\ell}(\bm{\Lambda}_{1}) = \diag\bigl(
            T_{\ell}(\eta(\lambda_{1}(\bm{\Lambda}_{1}))), \dots,
        T_{\ell}(\eta(\lambda_{d}(\bm{\Lambda}_{1}))) \bigr) \in \R^{d \times d}
    \]
    is positive definite and hence invertible. We therefore define
    \[
        \bm{\Delta}_{j}
        = \psi_{j}(\bm{\Lambda}_{1})
        (\psi_{j - 1}(\bm{\Lambda}_{1}))^{-1}
        = \diag\biggl(
            \frac{T_j(\eta(\lambda_1(\bm{\Lambda}_{1})))}{T_{j-1}(\eta(\lambda_1(\bm{\Lambda}_{1})))},
            \dots,
            \frac{T_j(\eta(\lambda_d(\bm{\Lambda}_{1})))}{T_{j-1}(\eta(\lambda_d(\bm{\Lambda}_{1})))}
        \biggr).
    \]
    Consider
    \[
        \mathbf{v}
        = \biggl(\mathbf{I} -
        \frac{\psi_{j}(\mathbf{K}^{2})}{\psi_{j}(0)} \biggr)
        \bm{\varepsilon}_{0}
        + \frac{\psi_{j - 1}(\mathbf{K}^{2})}{\psi_{j}(0)}
        \mathbf{F} \bm{\Delta}_{j} \bm{\xi}_{1}, \quad
        \psi_{j}(0) = T_{j} \Bigl(\frac{L + \tau_{2}}{L - \tau_{2}} \Bigr).
    \]
    As in the proof of \Cref{thm:full-spectrum}, the first term
    on the right-hand side belongs to
    $\range(\mathbf{W}_{2j}) = \Kry_{j}(\mathbf{K}^{2}, \mathbf{D})$.
    Since $\range(\mathbf{F}) \subseteq \range(\mathbf{D})$, the second
    term also belongs to $\range(\mathbf{W}_{2j})$.
    Hence $\mathbf{v} \in \range(\mathbf{W}_{2j})$.
    Using \eqref{eq:spectral-block-relations},
    $\mathbf{X}_{1}^{\top} \mathbf{F} = \mathbf{I}_{d}$, and
    $\mathbf{X}_{2}^{\top} \mathbf{F} = \bm{\Gamma}$, we obtain
    \begin{align*}
        \mathbf{X}_{1}^{\top}(\bm{\varepsilon}_{0} - \mathbf{v})
        & = \frac{1}{\psi_{j}(0)} \bigl(\psi_{j}(\bm{\Lambda}_{1})
        - \psi_{j - 1}(\bm{\Lambda}_{1}) \bm{\Delta}_{j} \bigr)
        \bm{\xi}_{1} = \mathbf{0},
        \shortintertext{and}
        \mathbf{X}_{2}^{\top}(\bm{\varepsilon}_{0} - \mathbf{v})
        & = \frac{1}{\psi_{j}(0)} \bigl(\psi_{j}(\bm{\Lambda}_{2}) \bm{\xi}_{2}
            - \psi_{j - 1}(\bm{\Lambda}_{2}) \bm{\Gamma}
        \bm{\Delta}_{j} \bm{\xi}_{1} \bigr).
    \end{align*}
    From \eqref{eq:initial-error-subspace} and
    \eqref{eq:krylov-spectral-subspace}, we have
    \[
        \bm{\varepsilon}_{0} - \mathbf{v}
        = \mathbf{X}_{1}
        \mathbf{X}_{1}^{\top}(\bm{\varepsilon}_{0} - \mathbf{v})
        + \mathbf{X}_{2} \mathbf{X}_{2}^{\top}(\bm{\varepsilon}_{0} -
        \mathbf{v})
        = \mathbf{X}_{2} \mathbf{X}_{2}^{\top}(\bm{\varepsilon}_{0} -
        \mathbf{v}).
    \]
    Hence, by \eqref{eq:spectral-block-relations},
    \begin{align}
        \lVert \bm{\varepsilon}_{0} - \mathbf{v} \rVert_{\mathbf{K}^{2}}
        & = \lVert \mathbf{X}_{2} \mathbf{X}_{2}^{\top}
        (\bm{\varepsilon}_{0} - \mathbf{v}) \rVert_{\mathbf{K}^{2}} \notag \\
        &= \lVert \bm{\Lambda}_{2}^{1/2} \mathbf{X}_{2}^{\top}
        (\bm{\varepsilon}_{0} - \mathbf{v}) \rVert \notag \\
        &= \frac{1}{\psi_{j}(0)} \bigl\lVert
        \bm{\Lambda}_{2}^{1/2} \psi_{j}(\bm{\Lambda}_{2}) \bm{\xi}_{2}
        - \bm{\Lambda}_{2}^{1/2} \psi_{j -
        1}(\bm{\Lambda}_{2}) \bm{\Gamma}
        \bm{\Delta}_{j} \bm{\xi}_{1} \bigr\rVert. \label{eq:X2err}
    \end{align}
    We next estimate \eqref{eq:X2err}. Since
    $\lambda(\bm{\Lambda}_{2}) \in [\tau_{2}, L]$, we have
    $\lambda(\eta(\bm{\Lambda}_{2})) \in [-1, 1]$.
    For every $\ell\ge0$, the inequality $|T_{\ell}(t)|\le1$ on $[-1,1]$ gives
    \[
        \lVert \psi_{\ell}(\boldsymbol{\Lambda}_{2}) \rVert
        \le\max_{t\in[\tau_2,L]}|T_{\ell}(\eta(t))|
        \le 1.
    \]
    Since $\bm{\Lambda}_{2}$ and $\psi_{\ell}(\bm{\Lambda}_{2})$ are both
    diagonal, for any vector $\mathbf{y}$, we have
    \begin{equation} \label{eq:LamPsi}
        \|\bm{\Lambda}_{2}^{1/2}
        \psi_{\ell}(\bm{\Lambda}_{2})\mathbf y\|
        =
        \|\psi_{\ell}(\bm{\Lambda}_{2})
        \bm{\Lambda}_{2}^{1/2}\mathbf y\|
        \le
        \|\bm{\Lambda}_{2}^{1/2}\mathbf y\|.
    \end{equation}
    By applying the triangle inequality and \eqref{eq:LamPsi} with
    $\ell=j$ and $\ell=j-1$ to \eqref{eq:X2err}, and then using
    the submultiplicativity of the norm, we obtain
    \begin{align}
        \lVert \bm{\varepsilon}_{0} - \mathbf{v} \rVert_{\mathbf{K}^{2}}
        &\le \frac{1}{\psi_j(0)}
        \bigl(
            \lVert \bm{\Lambda}_{2}^{1/2} \bm{\xi}_{2} \rVert
            + \lVert \bm{\Lambda}_{2}^{1/2}
            \bm{\Gamma}\bm{\Delta}_j
            \bm{\xi}_{1} \rVert
        \bigr) \notag \\
        &\le\frac{1}{\psi_j(0)}
        \bigl(
            \|\bm{\Lambda}_{2}^{1/2} \bm{\xi}_{2}\|
            +
            \|\bm{\Lambda}_{2}^{1/2} \bm{\Gamma}\|
            \,\|\bm{\Delta}_j \| \, \| \bm{\xi}_{1}\|
        \bigr). \label{eq:weighted-coupling-error}
    \end{align}
    Define
    \[
        \alpha_{i}=\eta(\lambda_{i}(\bm{\Lambda}_{1})),\quad
        \rho_{i}=\alpha_{i}+\sqrt{\alpha_{i}^{2}-1}.
    \]
    Since $\lambda_{i}(\bm{\Lambda}_{1})\ge1$ and $\eta$ is decreasing,
    \[
        1\le\alpha_{i}\le\eta(1)
        =\frac{L+\tau_{2}-2}{L-\tau_{2}}
        =\frac{\sigma_{1}^{2}+\sigma_{r-1}^{2}}
        {\sigma_{1}^{2}-\sigma_{r-1}^{2}}.
    \]
    The function $x+\sqrt{x^{2}-1}$ is increasing on $[1,\infty)$, so
    \[
        1\le\rho_{i}\le\eta(1)+\sqrt{\eta(1)^{2}-1}
        =\frac{\sigma_{1}+\sigma_{r-1}}{\sigma_{1}-\sigma_{r-1}}
        \eqcolon\rho.
    \]
    By \eqref{eq:chebyshev-explicit} and $\rho_{i}\ge1$, for each
    $1\le i\le d$,
    \[
        \frac{T_{j}(\alpha_{i})}{T_{j-1}(\alpha_{i})}
        =\frac{\rho_{i}^{j}+\rho_{i}^{-j}}
        {\rho_{i}^{j-1}+\rho_{i}^{-j+1}}
        =\rho_{i}\frac{1+\rho_{i}^{-2j}}
        {1+\rho_{i}^{2}\rho_{i}^{-2j}}
        \le\rho_{i}.
    \]
    The last inequality holds because $\rho_{i}\ge1$, so the denominator
    of the fraction is no smaller than its numerator. Therefore,
    \begin{equation}\label{eq:low-spectrum-ratio-bound}
        \|\bm{\Delta}_{j}\|
        =\max_{1\le i\le d}
        \frac{T_{j}(\alpha_{i})}{T_{j-1}(\alpha_{i})}
        \le\max_{1\le i\le d}\rho_{i}
        \le\rho.
    \end{equation}
    Combining \eqref{eq:weighted-coupling-error} and
    \eqref{eq:low-spectrum-ratio-bound} with the Chebyshev estimate
    \eqref{eq:chebbound}, with $\mu$ replaced by $\tau_{2}$, gives
    \begin{equation}\label{eq:weighted-coefficient-bound}
        \begin{aligned}
            \|\bm{\varepsilon}_{0}-\mathbf{v}\|_{\mathbf{K}^{2}}
            &\le T_{j} \Bigl(\frac{L + \tau_{2}}{L - \tau_{2}}
            \Bigr)^{-1} \bigl(
                \|\bm{\Lambda}_{2}^{1/2}\bm{\xi}_{2}\|
                +\rho\|\bm{\Lambda}_{2}^{1/2}\bm{\Gamma}\|
            \,\|\bm{\xi}_{1}\|\bigr) \\
            &\le 2\bigl(
                \|\bm{\Lambda}_{2}^{1/2}\bm{\xi}_{2}\|
                +\rho\|\bm{\Lambda}_{2}^{1/2}\bm{\Gamma}\|
            \,\|\bm{\xi}_{1}\|\bigr)
            \biggl(\frac{\sqrt{L/\tau_{2}}-1}
            {\sqrt{L/\tau_{2}}+1}\biggr)^{j}.
        \end{aligned}
    \end{equation}
    Since $\mathbf{v}\in\range(\mathbf{W}_{2j})$,
    \Cref{thm:transfer} and \eqref{eq:weighted-coefficient-bound}
    prove the stated bounds
    for the TriCG error and the TriMR residual norms.
\end{proof}

\section{Convergence bounds for odd iterations}
\label{sec:odd}
Since $\range(\mathbf{W}_{2j}) \subseteq
\range(\mathbf{W}_{2j + 1})$,
\Cref{thm:transfer} implies that the even-iteration bounds in
\Cref{sec:even} also bound the TriCG error and the TriMR residual
at the next odd iteration. However, these bounds do not account for
the larger subspace $\range(\mathbf{W}_{2j+1})$.

For a fixed $\widetilde{\mathbf{s}} \in \range(\mathbf{W}_{1})$,
the corresponding error and residual are
\begin{equation*}
    \widetilde{\bm{\varepsilon}} = \bm{\varepsilon}_{0} -
    \widetilde{\mathbf{s}}, \quad
    \widetilde{\mathbf{r}} = \mathbf{r}_{0} - \mathbf{K} \widetilde{\mathbf{s}}
    =
    \begin{bmatrix}
        \widetilde{\mathbf{r}}_{b} \\
        \widetilde{\mathbf{r}}_{c}
    \end{bmatrix}.
\end{equation*}
Thus, $\widetilde{\bm{\varepsilon}}$ is the exact solution of the
correction equation $\mathbf{K}\bm{\delta}=\widetilde{\mathbf{r}}$.
Using the residual blocks, define
\begin{equation*}
    \mathbf{D}_{\widetilde{\mathbf{s}}} =
    \begin{bmatrix}
        \widetilde{\mathbf{r}}_{b} & \mathbf{A} \widetilde{\mathbf{r}}_{c} &
        \mathbf{0} & \mathbf{0} \\
        \mathbf{0} & \mathbf{0} & - \widetilde{\mathbf{r}}_{c} &
        \mathbf{A}^{\top} \widetilde{\mathbf{r}}_{b}
    \end{bmatrix}.
\end{equation*}
By construction, $\mathbf{D}_{\widetilde{\mathbf{s}}} \one_{4}
= \mathbf{K} \widetilde{\mathbf{r}} = \mathbf{K}^{2}
\widetilde{\bm{\varepsilon}}$.
For $j\ge1$, the associated correction space is
$\Kry_{j}(\mathbf{K}^{2}, \mathbf{D}_{\widetilde{\mathbf{s}}})$.
By \eqref{eq:allspaces},
\[
    \range(\mathbf{K} \mathbf{W}_{k}) \subseteq \range(\mathbf{W}_{k
    + 1}), \quad
    \range(\mathbf{K}^{2} \mathbf{W}_{k}) \subseteq \range(\mathbf{W}_{k + 2}).
\]
Since $\mathbf{r}_{0}$ and $\widetilde{\mathbf{s}} \in \range(\mathbf{W}_{1})$,
we have $\widetilde{\mathbf{r}} \in \range(\mathbf{W}_{2})$. From
\[
    \begin{bmatrix}
        \mathbf{A} \widetilde{\mathbf{r}}_{c} & \mathbf{0} \\
        \mathbf{0} & \mathbf{A}^{\top} \widetilde{\mathbf{r}}_{b}
    \end{bmatrix} =
    \begin{bmatrix}
        \mathbf{0} & \mathbf{A} \\
        \mathbf{A}^{\top} & \mathbf{0}
    \end{bmatrix}
    \begin{bmatrix}
        \mathbf{0} & \widetilde{\mathbf{r}}_{b} \\
        \widetilde{\mathbf{r}}_{c} & \mathbf{0}
    \end{bmatrix},
\]
it follows that
\[
    \range\bigl(\mathbf{D}_{\widetilde{\mathbf{s}}}
        \begin{bmatrix} \mathbf{e}_{2} & \mathbf{e}_{4}
    \end{bmatrix}\bigr) \subseteq \range(\mathbf{K}_{0} \mathbf{W}_{2})
    \subseteq \range(\mathbf{W}_{3}).
\]
Hence $\range(\mathbf{D}_{\widetilde{\mathbf{s}}}) \subseteq
\range(\mathbf{W}_{3})$, and
\[
    \range \bigl((\mathbf{K}^{2})^{\ell}
    \mathbf{D}_{\widetilde{\mathbf{s}}} \bigr) \subseteq
    \range(\mathbf{W}_{2 \ell + 3}) \subseteq
    \range(\mathbf{W}_{2j + 1}), \quad 0 \le \ell \le j - 1.
\]
By the definition of the block Krylov subspace,
\begin{equation}
    \label{eq:cgoddaugmentation}
    \range(\mathbf{W}_{1}) + \Kry_{j}(\mathbf{K}^{2},
    \mathbf{D}_{\widetilde{\mathbf{s}}})
    \subseteq \range(\mathbf{W}_{2j + 1}).
\end{equation}

From \eqref{eq:commoncomparison} and \eqref{eq:cgoddaugmentation},
\begin{equation}
    \label{eq:cgoddtransfer}
    \begin{aligned}
        \lVert \bm{\varepsilon}_{2j + 1}^{\TCG} \rVert
        & \le \lVert \mathbf{r}_{2j + 1}^{\TMR} \rVert
        = \min_{\mathbf{v} \in \range(\mathbf{W}_{2j + 1})}
        \lVert \bm{\varepsilon}_{0} - \mathbf{v} \rVert_{\mathbf{K}^{2}} \\
        & \le \min_{\mathbf{v} \in \range(\mathbf{W}_{1}) +
        \Kry_{j}(\mathbf{K}^{2}, \mathbf{D}_{\widetilde{\mathbf{s}}})}
        \lVert \bm{\varepsilon}_{0} - \mathbf{v} \rVert_{\mathbf{K}^{2}} \\
        & \le \min_{\bm{\delta} \in
        \Kry_{j}(\mathbf{K}^{2}, \mathbf{D}_{\widetilde{\mathbf{s}}})}
        \lVert (\bm{\varepsilon}_{0} - \widetilde{\mathbf{s}}) -
        \bm{\delta} \rVert_{\mathbf{K}^{2}} \\
        & = \min_{\bm{\delta} \in
        \Kry_{j}(\mathbf{K}^{2}, \mathbf{D}_{\widetilde{\mathbf{s}}})}
        \lVert \widetilde{\bm{\varepsilon}} - \bm{\delta}
        \rVert_{\mathbf{K}^{2}}.
    \end{aligned}
\end{equation}
Using \eqref{eq:cgoddtransfer}, we obtain analogous convergence bounds
for odd iterations by following the analysis in \Cref{sec:even},
with $\Kry_{j}(\mathbf{K}^{2}, \mathbf{D}_{\widetilde{\mathbf{s}}})$,
$\mathbf{D}_{\widetilde{\mathbf{s}}}$, $\widetilde{\bm{\varepsilon}}$, and
$\widetilde{\mathbf{r}}$ playing the roles of
$\range(\mathbf{W}_{2j})$, $\mathbf{D}$, $\bm{\varepsilon}_{0}$, and
$\mathbf{r}_{0}$, respectively. Taking $\widetilde{\mathbf{s}} =
\mathbf{s}_{1}^{\TCG}$ or $\widetilde{\mathbf{s}} =
\mathbf{s}_{1}^{\TMR}$ gives bounds for the TriCG error or the TriMR
residual, respectively. The assumptions of the corresponding results
in \Cref{sec:even} must hold after these substitutions.

\section{Extension to general symmetric quasi-definite linear systems}
\label{sec:extension}
We now return to the general system \eqref{eq:general-system} and extend
the convergence bounds in \Cref{sec:even,sec:odd}
to this case. Let $\mathbf{H} = \diag(\mathbf{M}, \mathbf{N})$.
Define
\[
    \widetilde{\mathbf{K}}
    = \mathbf{H}^{-1/2} \mathbf{K} \mathbf{H}^{-1/2}
    =
    \begin{bmatrix}
        \mathbf{I}_{m} & \widetilde{\mathbf{A}} \\
        \widetilde{\mathbf{A}}^{\top} & - \mathbf{I}_{n}
    \end{bmatrix}, \quad
    \widetilde{\mathbf{A}} = \mathbf{M}^{-1/2} \mathbf{A}
    \mathbf{N}^{-1/2},
\]
and
\[
    \widetilde{\mathbf{s}} = \mathbf{H}^{1/2} \mathbf{s}, \quad
    \widetilde{\mathbf{f}} = \mathbf{H}^{-1/2} \mathbf{f}.
\]
The transformed system $\widetilde{\mathbf{K}} \widetilde{\mathbf{s}}
= \widetilde{\mathbf{f}}$ has the form \eqref{eq:identity-system}, with
starting vectors $\widetilde{\mathbf{b}} = \mathbf{M}^{-1/2} \mathbf{b}$
and $\widetilde{\mathbf{c}} = \mathbf{N}^{-1/2} \mathbf{c}$.
Applying \Cref{alg:orth_tri} to this system yields the
orthonormal bases $\widetilde{\mathbf{U}}_{k} = \mathbf{M}^{1/2} \mathbf{U}_{k}$
and $\widetilde{\mathbf{V}}_{k} = \mathbf{N}^{1/2} \mathbf{V}_{k}$.
Hence \eqref{eq:allspaces} yields
\begin{equation*}
    \widetilde{\mathbf{W}}_{k} = \mathbf{H}^{1/2} \mathbf{W}_{k}, \quad
    \range(\widetilde{\mathbf{W}}_{k})
    = \Kry_{k}(\widetilde{\mathbf{K}}, \widetilde{\mathbf{B}}_{0}), \quad
    \widetilde{\mathbf{B}}_{0} = \mathbf{H}^{-1/2} \mathbf{B}_{0}.
\end{equation*}
From the Galerkin condition \eqref{eq:tricg_galerkin} of TriCG and the
minimal residual property \eqref{eq:trimrsmall} of TriMR, we obtain
\[
    \mathbf{W}_{k}^{\top}(\mathbf{f} - \mathbf{K} \mathbf{s}_{k})
    = \widetilde{\mathbf{W}}_{k}^{\top}
    (\widetilde{\mathbf{f}} - \widetilde{\mathbf{K}} \mathbf{H}^{1/2}
    \mathbf{s}_{k}), \quad
    \lVert \mathbf{f} - \mathbf{K} \mathbf{s}_{k} \rVert_{\mathbf{H}^{-1}}
    = \lVert \widetilde{\mathbf{f}} - \widetilde{\mathbf{K}}
    \mathbf{H}^{1/2} \mathbf{s}_{k} \rVert.
\]
Thus, for either method, the iterates for the transformed system satisfy
$\widetilde{\mathbf{s}}_{k} = \mathbf{H}^{1/2} \mathbf{s}_{k}$. The error
and residual satisfy
\begin{equation}
    \label{eq:transformed-norm-identities}
    \widetilde{\bm{\varepsilon}}_{k} = \mathbf{H}^{1/2}
    \bm{\varepsilon}_{k}, \quad
    \widetilde{\mathbf{r}}_{k} = \mathbf{H}^{-1/2} \mathbf{r}_{k}, \quad
    \lVert \widetilde{\bm{\varepsilon}}_{k} \rVert = \lVert
    \bm{\varepsilon}_{k} \rVert_{\mathbf{H}}, \quad
    \lVert \widetilde{\mathbf{r}}_{k} \rVert = \lVert
    \mathbf{r}_{k} \rVert_{\mathbf{H}^{-1}}.
\end{equation}
Define the symmetric positive definite matrix $\mathbf{G} = \mathbf{K}
\mathbf{H}^{-1} \mathbf{K}$. By \Cref{thm:transfer},
\begin{equation*}
    \begin{aligned}
        \lVert \bm{\varepsilon}_{k}^{\TCG} \rVert_{\mathbf{H}}
        & = \lVert \widetilde{\bm{\varepsilon}}_{k}^{\TCG} \rVert \\
        & \le \min_{\widetilde{\mathbf{v}} \in
        \range(\widetilde{\mathbf{W}}_{k})}
        \lVert \widetilde{\bm{\varepsilon}}_{0} - \widetilde{\mathbf{v}}
        \rVert_{\widetilde{\mathbf{K}}^{2}} \\
        & = \min_{\mathbf{v} \in \range(\mathbf{W}_{k})}
        \lVert \mathbf{H}^{1/2}(\bm{\varepsilon}_{0} - \mathbf{v})
        \rVert_{\widetilde{\mathbf{K}}^{2}} \\
        & = \min_{\mathbf{v} \in \range(\mathbf{W}_{k})}
        \lVert \bm{\varepsilon}_{0} - \mathbf{v} \rVert_{\mathbf{G}} \\
        & = \lVert \widetilde{\mathbf{r}}_{k}^{\TMR} \rVert
        = \lVert \mathbf{r}_{k}^{\TMR} \rVert_{\mathbf{H}^{-1}}.
    \end{aligned}
\end{equation*}

These relations extend the convergence analysis in
\Cref{sec:even,sec:odd} to the general system
\eqref{eq:general-system} in the corresponding weighted norms.
For example, let $\widetilde{\mu}$ and
$\widetilde{L}$ be the smallest and largest eigenvalues of
$\widetilde{\mathbf{K}}^{2}$. Applying
\eqref{eq:cgunconditional} and \eqref{eq:mrunconditional} to the transformed
system and using \eqref{eq:transformed-norm-identities}, we obtain
\begin{align*}
    \lVert \bm{\varepsilon}_{2j}^{\TCG} \rVert_{\mathbf{H}}
    & \le 2 \lVert \bm{\varepsilon}_{0} \rVert_{\mathbf{G}}
    \left(\frac{\sqrt{\widetilde{L}/\widetilde{\mu}} - 1}
    {\sqrt{\widetilde{L}/\widetilde{\mu}} + 1} \right)^{j},
    \shortintertext{and}
    \lVert \mathbf{r}_{2j}^{\TMR} \rVert_{\mathbf{H}^{-1}}
    & \le 2 \lVert \mathbf{r}_{0} \rVert_{\mathbf{H}^{-1}}
    \left(\frac{\sqrt{\widetilde{L}/\widetilde{\mu}} - 1}
    {\sqrt{\widetilde{L}/\widetilde{\mu}} + 1} \right)^{j}.
\end{align*}
Similarly, the bounds \eqref{eq:bound2} and \eqref{eq:conditionalCG},
together with their odd-iteration counterparts in \Cref{sec:odd},
extend to the general system \eqref{eq:general-system}.

\section{Numerical experiments}
\label{sec:numerical}

In this section, we conduct numerical experiments to verify the
convergence bounds \eqref{eq:bound2} and \eqref{eq:conditionalCG} for the
TriCG error and TriMR residual norms, together with their odd-iteration
counterparts derived in \Cref{sec:odd}.
All computations use double precision and a zero initial approximation.
The basis vectors are generated with full reorthogonalization to limit
the loss of orthogonality due to round-off errors. Both methods terminate
when the residual norm falls below $10^{-12}$.

For the first two examples, we set $\mathbf{M} = \mathbf{I}$ and
$\mathbf{N} = \mathbf{I}$. We first consider the $200 \times 200$
diagonal matrix
\[
    \mathbf{A}=\diag(0,t_{1},\ldots,t_{199}),
    \quad t_{i}=2+\frac{8(i-1)}{198}.
\]
The entries of $\mathbf{b}$ and $\mathbf{c}$ are drawn independently
from the standard normal distribution and are then normalized so that
$\|[\mathbf{b}^{\top},\mathbf{c}^{\top}]^{\top}\|=1$.
The two smallest positive singular values are
\[
    \sigma_{r}=2,\quad
    \sigma_{r-1}=2+\frac{8}{198}\approx2.04.
\]

\begin{figure}[htb]
    \centering
    \includegraphics[width=0.92\textwidth]{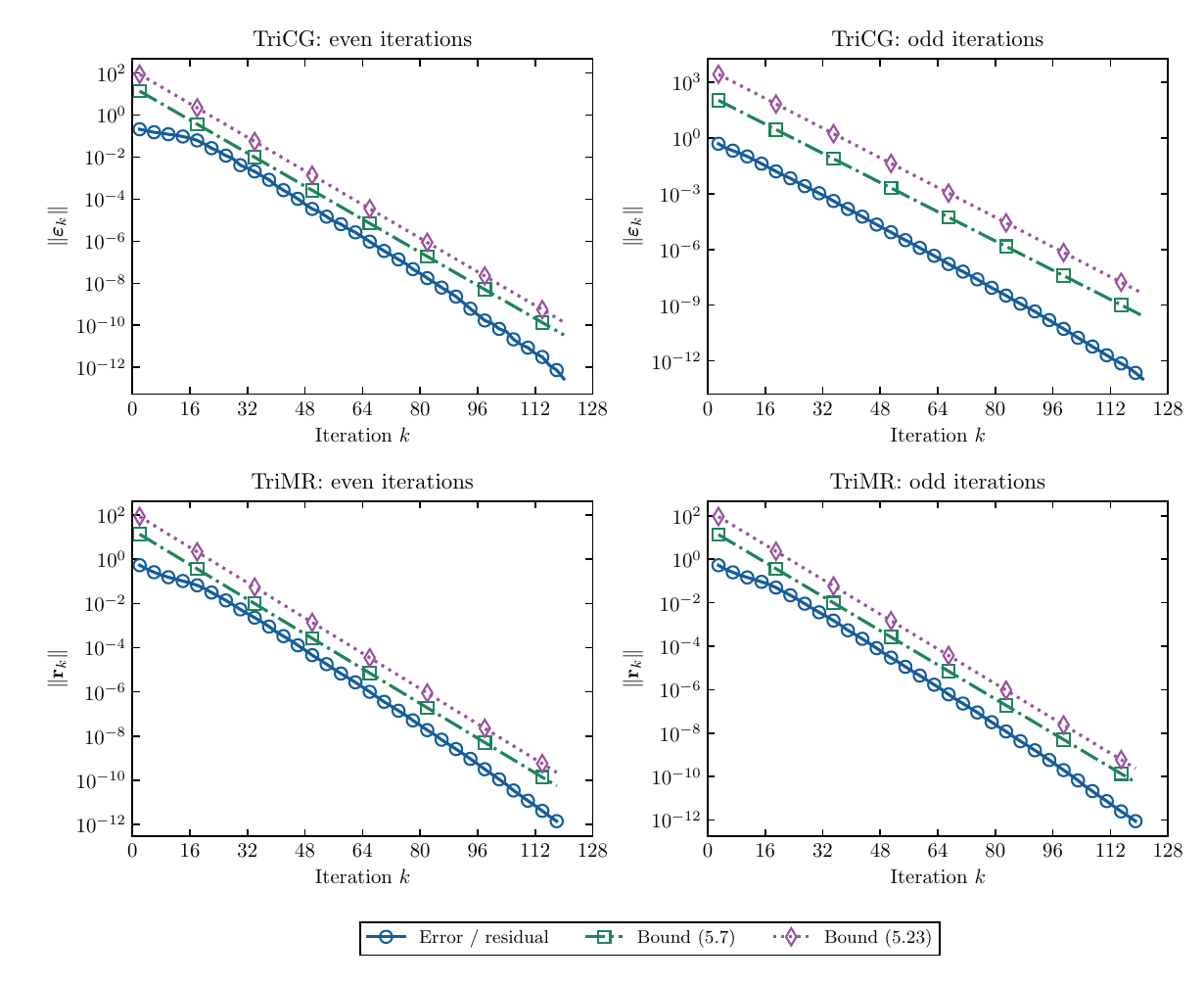}
    \caption{Convergence histories of TriCG and TriMR, with
        bounds~\eqref{eq:bound2} and \eqref{eq:conditionalCG} at even iterations
    and their counterparts from \Cref{sec:odd} at odd iterations.}
    \label{fig:exp-uniform}
\end{figure}

\Cref{fig:exp-uniform} displays the convergence histories of TriCG and TriMR
at even and odd iterations, together with the corresponding convergence
bounds. The convergence factors of bounds~\eqref{eq:bound2} and
\eqref{eq:conditionalCG} are approximately $0.64$ and $0.63$,
respectively, and are therefore very close. The figure shows that
the curves of error and residual norms are approximately parallel to
those of the two bounds. For even iterations, the corresponding
constants are $C\approx21.70$ and $B\approx137.92$. Since the convergence
factors are very close, the smaller constant makes
bound~\eqref{eq:bound2} tighter than bound~\eqref{eq:conditionalCG}
over the displayed even iterations. Similar behavior is observed at
odd iterations.

\FloatBarrier

We next consider the $200\times200$ diagonal matrix
\[
    \mathbf{A}=\diag(0,2,t_{1},\ldots,t_{198}),
    \quad t_{i}=6+\frac{4(i-1)}{197}.
\]
We use the same vectors $\mathbf{b}$ and $\mathbf{c}$ as in the first example.
The two smallest positive singular values are
\[
    \sigma_{r}=2,\quad \sigma_{r-1}=6,
\]
giving a larger gap than in the first example.

\begin{figure}[htb]
    \centering
    \includegraphics[width=0.92\textwidth]{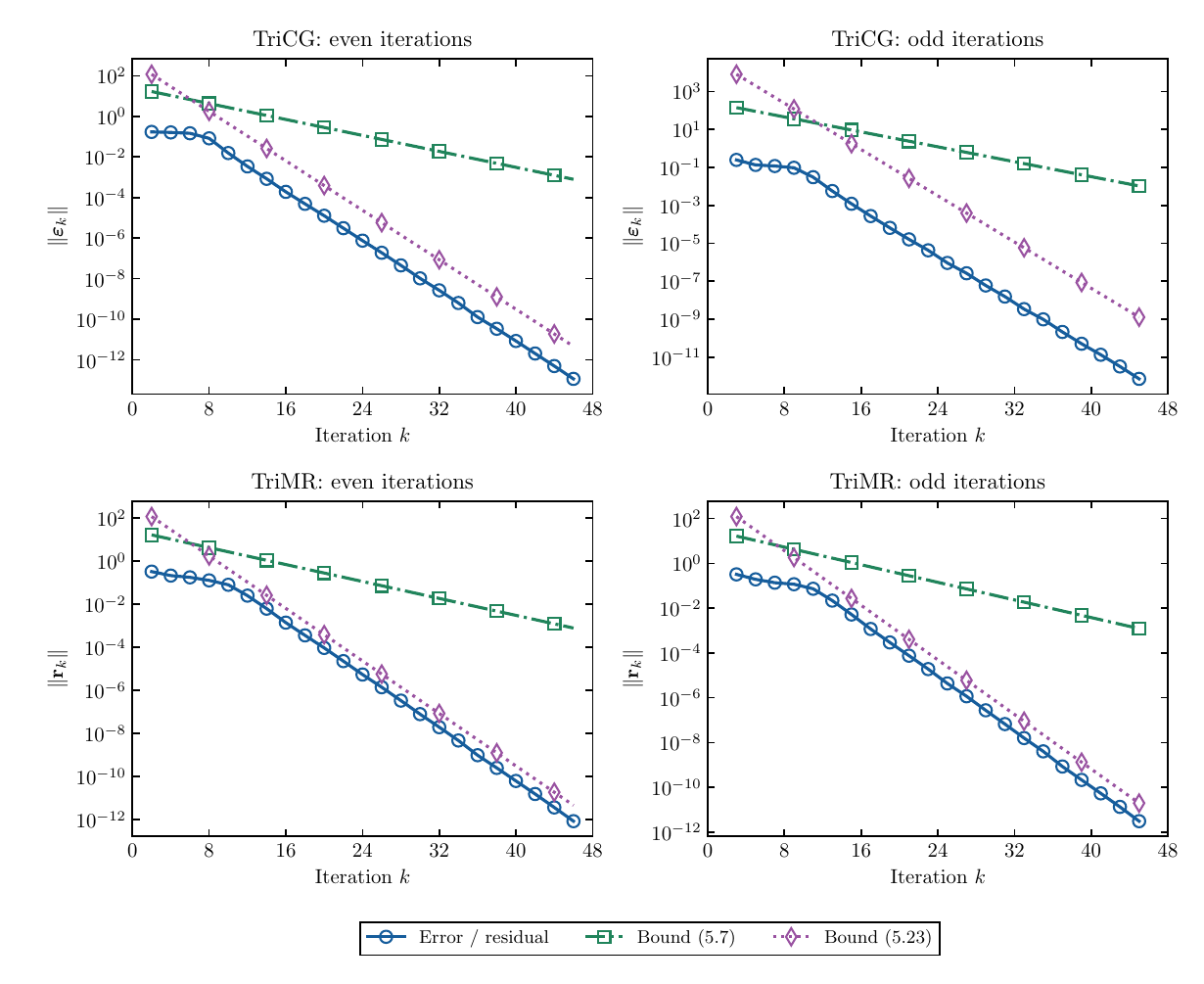}
    \caption{Convergence histories of TriCG and TriMR, with
        bounds~\eqref{eq:bound2} and \eqref{eq:conditionalCG} at even iterations
    and their counterparts from \Cref{sec:odd} at odd iterations.}
    \label{fig:exp-spectrum}
\end{figure}

\Cref{fig:exp-spectrum} compares the even- and odd-iteration convergence
histories of TriCG and TriMR with the corresponding bounds.
The convergence factors of bounds~\eqref{eq:bound2} and
\eqref{eq:conditionalCG} are approximately $0.64$ and $0.25$,
respectively, with the latter being much smaller.
The figure shows that both bounds remain above the computed norms.
After the first few iterations, bound~\eqref{eq:conditionalCG} follows
the decay of the computed norms more closely, whereas
bound~\eqref{eq:bound2} decreases much
more slowly. For even iterations, the constants are $C\approx26.15$
and $B\approx480.23$. Although bound~\eqref{eq:conditionalCG} is initially
larger because of its larger constant, its faster decay makes it tighter
than bound~\eqref{eq:bound2} from $k=8$ onward.
At odd iterations, bound~\eqref{eq:conditionalCG} likewise becomes tighter
than bound~\eqref{eq:bound2} as the iteration proceeds.

Finally, we consider the general SQD linear system arising at the initial
interior-point iterate for the convex quadratic optimization
problem \texttt{mosarqp2}, taken from the SQD
collection \cite{orban2015acollection}.
Here $\mathbf{A}\in\mathbb{R}^{1500\times2400}$, $\mathbf{M}=\mathbf{I}$,
and $\mathbf{N}\in\mathbb{R}^{2400\times2400}$ is symmetric positive definite
with $\mathbf{N}\ne\mathbf{I}$. The two smallest singular values of
$\mathbf{M}^{-1/2}\mathbf{A}\mathbf{N}^{-1/2}$ are approximately
$0.27$ and $0.29$, respectively.

\begin{figure}[htb]
    \centering
    \includegraphics[width=0.92\textwidth]{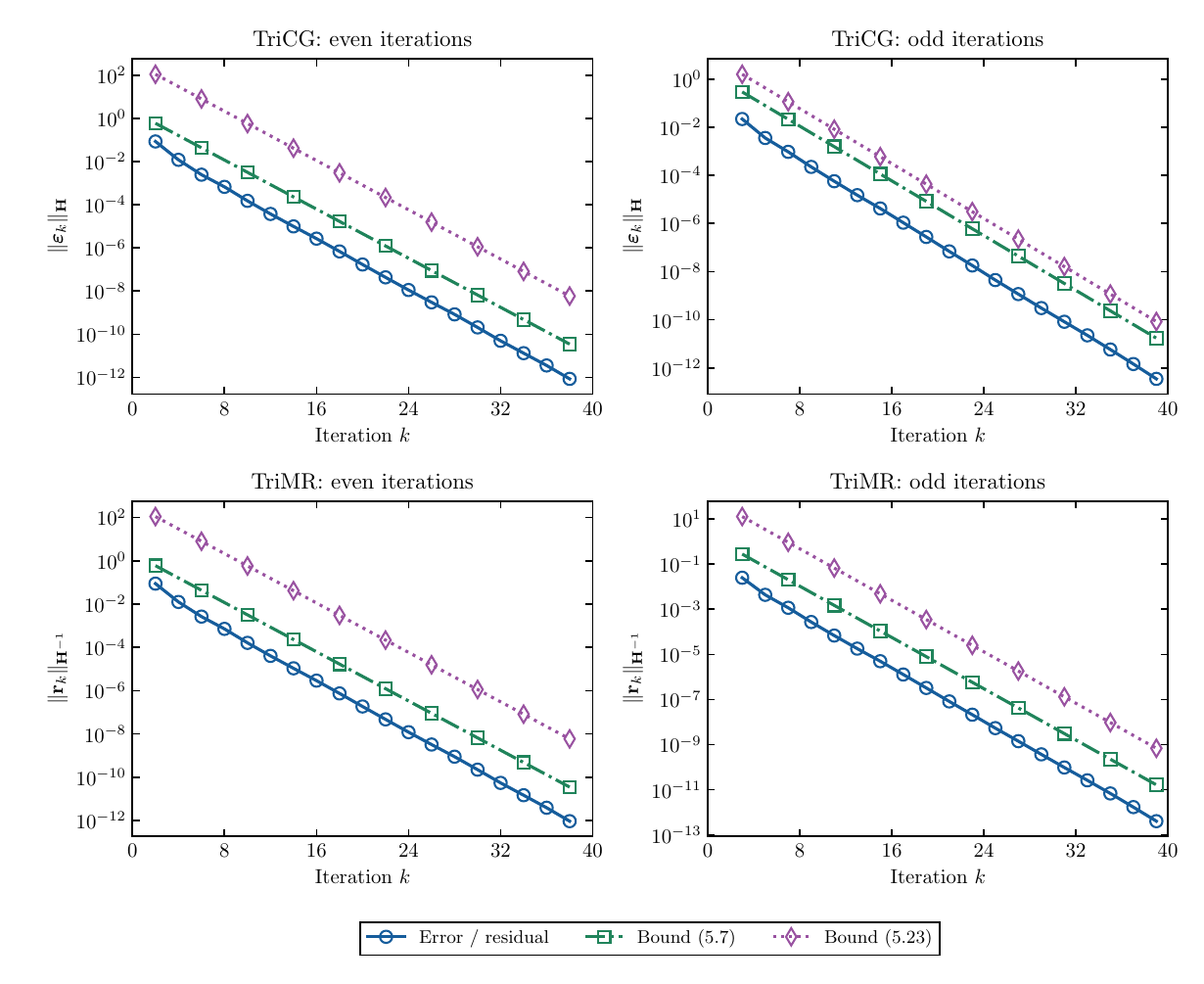}
    \caption{Convergence histories of TriCG and TriMR for \texttt{mosarqp2},
        with the weighted forms of bounds~\eqref{eq:bound2} and
        \eqref{eq:conditionalCG} at even iterations and their counterparts
    at odd iterations, as described in \Cref{sec:extension}.}
    \label{fig:exp-sqd-0}
\end{figure}

\Cref{fig:exp-sqd-0} presents the convergence histories of the TriCG
error norm and the TriMR residual norm at even and odd iterations,
along with the corresponding bounds. The convergence factors of
bounds~\eqref{eq:bound2} and \eqref{eq:conditionalCG} are very close,
both being approximately $0.27$. The figure shows that both
bounds lie above the computed error and residual norms and decay at
rates similar to those of the computed norms.
For even iterations, the corresponding constants are $C\approx2.18$
and $B\approx407.40$. Since the convergence factors are very close, the smaller
constant makes bound~\eqref{eq:bound2} tighter than
bound~\eqref{eq:conditionalCG} in \Cref{fig:exp-sqd-0}.
Similar behavior is observed at odd iterations.

\section{Concluding remarks and future work}
\label{sec:discussion}
In this paper, we present the convergence analysis of TriCG and TriMR
for SQD linear systems.
For the special case~\eqref{eq:identity-system}, we establish that the
TriCG error norm is bounded by the TriMR residual norm and use this
relation to derive convergence bounds for both methods. We then extend
these bounds to the general system~\eqref{eq:general-system} in the
corresponding weighted norms.

For even iterations, we derive a CG-type bound and two refined bounds
that exploit the block structure of the system. The first refinement
uses the smallest and largest positive singular values of the
off-diagonal block to determine the convergence factor. Under
additional spectral and rank conditions,
the second refinement replaces the smallest positive singular value in
the convergence factor with the second-smallest, potentially reducing
the factor further. For odd iterations, we obtain corresponding
bounds by applying the same analysis to correction problems formed from
the first TriCG or TriMR iterate. Numerical experiments support our
theoretical results.

Our analysis assumes exact arithmetic and exact solution of the linear
systems involving the diagonal blocks. A direction for future work is
to exploit the block structure to analyze the effects of
round-off errors and inexact solves on the convergence of TriCG and TriMR.
Related analyses of finite-precision Lanczos and CG methods and inexact
Krylov methods can be found in \cite{meurant2006lanczos} and
\cite{simoncini2003theory}, respectively.

\bibliographystyle{abbrv}
\begingroup
\small
\bibliography{conv}
\endgroup
\end{document}